\documentclass{article}
\usepackage{amssymb}
\usepackage{graphicx}
\usepackage{amsmath}

\input{tcilatex}
\begin{document}

\begin{center}
\textbf{Operator-valued multipliers in vector-valued weighted Besov spaces
and applications}
\end{center}

\QTP{Body Math}
\ 

\begin{center}
{\textbf{Veli B. Shakhmurov}}

Okan University, Department of Mechanical engineering, Akfirat, Tuzla 34959
Istanbul, Turkey,

E-mail: veli.sahmurov@okan.edu.tr

\textbf{Rishad Shahmurov}

\bigskip E-mail: shahmurov@hotmail.com\ \ \ \ \ \ \ \ \ \ \ \ \ \ 
\end{center}

\QTP{Body Math}
\ \ \ \ \ \ \ \ \ \ \ \ \ \ \ \ \ \ \ \ \ \ \ \ \ \ \ \ \ \ \ \ \ \ \ \ \ \
\ \ \ \ \ \ \ \ \ \ \ \ \ \ 

\begin{center}
\textbf{ABSTRACT}
\end{center}

\begin{quote}
\ \ \ \ \ \ \ \ \ \ \ \ \ \ \ 
\end{quote}

\ \ \ The operator-valued multiplier theorems in $E-$valued Besov spaces are
studied, where $E_{0}$, $E$ are two Banach spaces and $E_{0}\subset E$.
These results permit us to show embedding theorems in $E_{0}$-valued
weighted Besov-Lions type spaces $B_{p,q,\gamma }^{l,s}\left( \Omega
;E_{0},E\right) .$ The most regular class of interpolation space $E_{\alpha
},$ between $E_{0}$ and $E$ are found such that the mixed differential
operator $D^{\alpha }$ is bounded from $B_{p,q,\gamma }^{l,s}\left( \Omega
;E_{0},E\right) $ to $B_{p,q,\gamma }^{s}\left( \Omega ;E_{\alpha }\right) $
and Ehrling-Nirenberg-Gagliardo type sharp estimates are established. By
using these results the separability properties of degenerate differential
operators are studied. Especially, we prove that the associated differential
operators are positive and also are generators of analytic semigroups.
Moreover, maximal $B_{p,q,\gamma }^{s}$-regularity properties for abstract
elliptic equation, Cauchy problem for degenerate abstract parabolic equation
and the infinite systems of degenerate parabolic equations are studied.\ \ \
\ 

\begin{center}
\bigskip\ \ \textbf{AMS:34G10, 35J25, 35J70}
\end{center}


\begin{center}
\textbf{1. Introduction }
\end{center}

Fourier multipliers in vector-valued function spaces has been studied e.g.
in $\left[ 17\right] ,$ $\left[ 28\right] ,$ $\left[ 32\right] .$
Operator-valued Fourier multipliers in weighted spaces have been
investigated in $\left[ 1\right] ,$ $\left[ 9-12\right] ,$ $\left[ 13\right]
,\left[ 30\right] .$ Mikhlin type Fourier multiplierers in scalar weighted
spaces have been studied e.g. in $\left[ 14\right] $ and $\left[ 30\right] $%
. Moreover, operator-valued Fourier multiplers in weighted abstract $L_{p}$
spaces were investigated e.g. in $\left[ 2\right] $, $\left[ 7\right] ,$ $%
\left[ 13\right] ,$ $\left[ 16\right] .$ Regularity properties of
differential-operator equations (DOEs) have been studied e.g. in $\left[ 1%
\right] $, $\left[ 3\right] ,$ $\left[ 9\right] ,$ $\left[ 21-26\right] $, $%
\left[ 30-31\right] .$ A comprehensive introduction to DOEs and historical
references may be found in $\left[ 1\right] $ and $\left[ 31\right] .$

In this paper, operator-valued multiplier theorems in $E-$valued weighted
Lebesque and Besov spaces are shown. Then we consider the $E$- valued
anisotropic Besov spaces $B_{p,q,\gamma }^{l,s}\left( \Omega ;E_{0},E\right) 
$, here $E_{0}$, $E$ are two Banach spaces, $E_{0}$ is continuously and
densely embedded into $E,$ and $\gamma =\gamma \left( x\right) $ is weighted
function from $A_{p}$, $p\in \left( 1,\infty \right) $ class. We prove
boundedness and compactness of embedding operators in these spaces. This
result generalized and improved the results $\left[ \text{4, \S\ 9, 27, \S\ %
1.7}\right] $ for scalar Sobolev space, the result $\left[ 15\right] $ for
one dimensional Sobolev-Lions spaces and the results $\left[ 22-23\right] $
for Hilbert-space valued class. Finally, we consider differential-operator
equation%
\begin{equation}
\ Lu=\sum\limits_{\left\vert \alpha \right\vert =2l}a_{\alpha }D^{\alpha
}u+Au+\sum\limits_{\left\vert \alpha \right\vert <2l}A_{\alpha }D^{\alpha
}u=f  \tag{1.1}
\end{equation}%
where $a_{\alpha }$ are complex numbers, $A$ and $A_{\alpha }\left( x\right) 
$ are linear operators in a Banach space $E$, $\alpha =\left( \alpha
_{1},\alpha _{2},...,\alpha _{n}\right) .$

We say that the problem $\left( 1.1\right) $ is $B_{p,q,\gamma }^{s}\left(
R^{n};E\right) $-separable, if there exists a unique solution 
\begin{equation*}
u\in B_{p,q,\gamma }^{2l,s}\left( \Omega ;E\left( A\right) ,E\right)
\end{equation*}%
of $\left( 1.1\right) $ for all $f\in B_{p,q,\gamma }^{s}\left(
R^{n};E\right) $ and there exists a positive constant $C$ independent of $f$
such that the coercive estimate holds 
\begin{equation}
\sum\limits_{\left\vert \alpha \right\vert =2l}\left\Vert D^{\alpha
}u\right\Vert _{B_{p,q,\gamma }^{s}\left( R^{n};E\right) }+\left\Vert
Au\right\Vert _{B_{p,q,\gamma }^{s}\left( R^{n};E\right) }\leq C\left\Vert
f\right\Vert _{B_{p,q,\gamma }^{s}\left( R^{n};E\right) }.  \tag{1.2}
\end{equation}

The estimate $\left( 1.2\right) $ implies that if $f\in B_{p,q,\gamma
}^{s}\left( R^{n};E\right) $ and $u$ is the solution of the problem $\left(
1.1\right) $ then all terms of the equation $\left( 1.1\right) $ belong to $%
B_{p,q,\gamma }^{s}\left( R^{n};E\right) $ (i.e. all terms are separable in $%
B_{p,q,\gamma }^{s}\left( R^{n};E\right) $)$.$

The above estimate implies that the inverse of the differential operator
generated by $\left( 1.1\right) $ is bounded from $B_{p,q,\gamma }^{s}\left(
R^{n};E\right) $ to%
\begin{equation*}
B_{p,q,\gamma }^{2l,s}\left( \Omega ;E\left( A\right) ,E\right) .
\end{equation*}

By using the separability properties of $\left( 1.1\right) $ we show the
maximal regularity properties of the following abstract parabolic Cauchy
problem%
\begin{equation}
\ \partial _{t}u+\sum\limits_{\left\vert \alpha \right\vert =2l}a_{\alpha
}D^{\alpha }u+Au=f\left( t,x\right) ,  \tag{1.3}
\end{equation}%
\begin{equation*}
u\left( 0,x\right) =0
\end{equation*}%
in weighted Besov spaces.

The paper is organized as follows. In Section 2 the necessary tools from
Banach space theory and some background materials are given. In Sections 3-5
the multiplier theorems in vector-valued weighted Lebesque and Besov spaces
are proved. In Sections 6-8 by using these multiplier theorems, embedding
theorems in $E$-valued weighted Besov type spaces are shown. Finally, in
Sections 9-14 the separability properties of problems $\left( 1.1\right) $, $%
\left( 1.3\right) $ and their applications are established.

\begin{center}
\textbf{2. Notations and background }
\end{center}

\ \ \ Let $E$ be a Banach space and $\gamma =\gamma \left( x\right) ,$ $%
x=\left( x_{1},x_{2},...,x_{n}\right) $ be a positive measurable function on
the measurable subset $\Omega \subset R^{n}.$ Let $L_{p,\gamma }\left(
\Omega ;E\right) $ denote the space of strongly measurable $E-$valued
functions that are defined on $\Omega $ with the norm

\begin{equation*}
\left\Vert f\right\Vert _{L_{p,\gamma }}=\left\Vert f\right\Vert
_{L_{p,\gamma }\left( \Omega ;E\right) }=\left( \int \left\Vert f\left(
x\right) \right\Vert _{E}^{p}\gamma \left( x\right) dx\right) ^{\frac{1}{p}},%
\text{ }1\leq p<\infty ,
\end{equation*}

\begin{equation*}
\left\Vert f\right\Vert _{L_{\infty ,\gamma }\left( \Omega ;E\right) }=\text{%
ess}\sup\limits_{x\in \Omega }\left\Vert f\left( x\right) \right\Vert
_{E}\gamma \left( x\right) ,\text{ }p=\infty .
\end{equation*}

For $\gamma \left( x\right) \equiv 1,$ the space $L_{p,\gamma }\left( \Omega
;E\right) $ will be denoted by $L_{p}=L_{p}\left( \Omega ;E\right) .$

The weight $\gamma $ is said to be satisfy an $A_{p}$ condition $\left[ 18%
\right] $, i.e., $\ \gamma \in A_{p},$ $1<p<\infty $ if there is a positive
constant $C$ such that 
\begin{equation*}
\left( \frac{1}{\left\vert Q\right\vert }\int\limits_{Q}\gamma \left(
x\right) dx\right) \left( \frac{1}{\left\vert Q\right\vert }%
\int\limits_{Q}\gamma ^{-\frac{1}{p-1}}\left( x\right) dx\right) ^{p-1}\leq
C,
\end{equation*}

for all cubes $Q\subset R^{n}.$

The Banach space\ $E$\ is called a UMD-space and written as $E\in $ UMD if
only if the Hilbert operator 
\begin{equation*}
\left( Hf\right) \left( x\right) =\lim\limits_{\varepsilon \rightarrow
0}\int\limits_{\left\vert x-y\right\vert >\varepsilon }\frac{f\left(
y\right) }{x-y}dy
\end{equation*}%
is bounded in the space $L_{p}\left( R,E\right) ,$ $p\in \left( 1,\infty
\right) $ (see e.g. $\left[ 6\right] $). UMD spaces include e.g. $L_{p}$, $%
l_{p}$ spaces and Lorentz spaces $L_{pq},$ $p,$ $q\in \left( 1,\infty
\right) $.

Let $\mathbb{C}$ be a set of complex numbers and\ 

\begin{equation*}
\ S_{\varphi }=\left\{ \xi ;\text{ \ }\xi \in \mathbb{C}\text{, \ }%
\left\vert \arg \xi \right\vert \leq \varphi \right\} \cup \left\{ 0\right\}
,0\leq \varphi <\pi .
\end{equation*}%
Let $E_{1}$ and\ $E_{2}$ be two Banach spaces. $B\left( E_{1},E_{2}\right) $
denotes the space of bounded linear operators ifrom $E_{1}$ to $E_{2}.$ For $%
E_{1}=E_{2}=E$ it will denote by $B\left( E\right) .$

A linear operator\ $A$ is said to be positive in\ a Banach\ space $E$,\ with
bound $M$ if\ $D\left( A\right) $ is dense on $E$ and 
\begin{equation*}
\left\Vert \left( A+\xi I\right) ^{-1}\right\Vert _{B\left( E\right) }\leq
M\left( 1+\left\vert \xi \right\vert \right) ^{-1}
\end{equation*}%
with $\xi \in S_{\varphi },\varphi \in \left[ 0,\right. \left. \pi \right) ,$
where $M$ is a positive constant and $I$ is an identity operator in $E.$
Sometimes instead of $A+\xi I$\ will be written $A+\xi $ and denoted by $%
A_{\xi }.$ It is known $\left[ \text{29, \S 1.15.1}\right] $ there exist
fractional powers\ $A^{\theta }$of the positive operator $A.$

\textbf{Definition 2.1}. A positive operator $A$ is said to be $R-$positive
in the Banach space $E$ if there exists $\varphi \in \left[ 0,\right. \left.
\pi \right) $ such that the set 
\begin{equation*}
\left\{ \left( \xi \right) \left( A+\xi I\right) ^{-1}:\xi \in S_{\varphi
}\right\}
\end{equation*}%
is $R$-bounded (see e.g. $\left[ 30\right] $).

$\sigma _{\infty }\left( E\right) $ will denote the space of compact
operators in $E.$

Let $E\left( A^{\theta }\right) $ denote the space $D\left( A^{\theta
}\right) $ with graphical norm defined as 
\begin{equation*}
\left\Vert u\right\Vert _{E\left( A^{\theta }\right) }=\left( \left\Vert
u\right\Vert ^{p}+\left\Vert A^{\theta }u\right\Vert ^{p}\right) ^{\frac{1}{p%
}},1\leq p<\infty ,-\infty <\theta <\infty .
\end{equation*}%
By $\left( E_{1},E_{2}\right) _{\theta ,p}$ will be denoted an interpolation
space\ obtained from $\left\{ E_{1},E_{2}\right\} $ by the $K-$method $\left[
\text{29, \S 1.3.1}\right] $, where $\theta \in \left( 0,1\right) ,$ $p\in %
\left[ 0,1\right] $. We denote by $D\left( R^{n};E\right) $ the space of $E-$%
valued $C^{\infty }-$ function with compact support, equipped with the usual
inductive limit topology and $S\left( E\right) =S\left( R^{n};E\right) $
denote the $E-$valued Schwartz space of rapidly decreasing, smooth
functions. For $E=\mathbb{C}$ we simply write $D\left( R^{n}\right) $ and $%
S=S\left( R^{n}\right) $, respectively. $D^{\prime }\left( R^{n};E\right) $ $%
=L\left( D\left( R^{n}\right) ,E\right) $ denote the space of $E-$valued
distributions and $S^{\prime }\left( E\right) =S^{\prime }\left(
R^{n};E\right) $ is a space of linear continued mapping from $S\left(
R^{n}\right) $ into\ $E.$ The Fourier transform for $u\in S^{\prime }\left(
R^{n};E\right) $ is defined by 
\begin{equation*}
F\left( u\right) \left( \varphi \right) =u\left( F\left( \varphi \right)
\right) \text{, }\varphi \in S\left( R^{n}\right) .
\end{equation*}%
Let $\gamma $ be such that $S\left( R^{n};E_{1}\right) $ is dense in $%
L_{p,\gamma }\left( R^{n};E_{1}\right) .$ A function 
\begin{equation*}
\Psi \in C^{\left( l\right) }\left( R^{n};B\left( E_{1},E_{2}\right) \right)
\end{equation*}
is called a multiplier from\ $L_{p,\gamma }\left( R^{n};E_{1}\right) $ to $%
L_{q,\gamma }\left( R^{n};E_{2}\right) $ if there exists a positive constant 
$C$ such that 
\begin{equation*}
\left\Vert F^{-1}\Psi \left( \xi \right) Fu\right\Vert _{L_{q,\gamma }\left(
R^{n};E_{2}\right) }\leq C\left\Vert u\right\Vert _{L_{p,\gamma }\left(
R^{n};E_{1}\right) }
\end{equation*}%
for all $u\in S\left( R^{n};E_{1}\right) $.

In a similar we can define the multiplier from $B_{p,q,\gamma }^{s}\left(
R^{n};E_{1}\right) $ to $B_{p,q,\gamma }^{s}\left( R^{n};E_{2}\right) .$

We denote the set of all multipliers fom\ $L_{p,\gamma }\left(
R^{n};E_{1}\right) $ to $L_{q,\gamma }\left( R^{n};E_{2}\right) $ by $%
M_{p,\gamma }^{q,\gamma }\left( E_{1},E_{2}\right) .$ For $E_{1}=E_{2}=E$ we
denote the $M_{p,\gamma }^{q,\gamma }\left( E_{1},E_{2}\right) $ by $%
M_{p,\gamma }^{q,\gamma }\left( E\right) .$

\textbf{Definition 2.2. }Let $\gamma $ be a positive measurable function on $%
R^{n}$. Assume $E$ is a Banach space and $p\in \left[ 1,2\right] .$ Suppose
there exists a positive constant $C_{0}=C_{0}\left( p,\gamma ,E\right) $ so
that 
\begin{equation}
\left\Vert Fu\right\Vert _{L_{p^{\prime },\gamma }\left( R^{n};E\right)
}\leq C_{0}\left\Vert Fu\right\Vert _{L_{p,\gamma }\left( R^{n};E\right) } 
\tag{2.1}
\end{equation}%
for $\frac{1}{p}+\frac{1}{p^{\prime }}=1$ and each $u\in S\left(
R^{n};E\right) .$ Then $E$ is called weighted Fourier type $\gamma $ and $p.$
It is called Fourier type $p\in \left[ 1,2\right] $ if $\gamma \left(
x\right) \equiv 1.$

\textbf{Remark 2.1.} The estimate $\left( 2.1\right) $ shows that each
Banach space $E$ has weighted Fourier type $\gamma $ and $1.$ By Bourgain $%
\left[ 6\right] $ has shown that each $B-$convex Banach space (thus, in
particular, each uniformly convex Banach space) has some non-trivial Fourier
type $p\in \left[ 1,2\right] $, i.e. $UMD$ spaces are Fourier type for some $%
p\in \left[ 1,2\right] .$

\bigskip In order to define abstract Besov spaces we consider the
dyadic-like subsets $\left\{ J_{k}\right\} _{k=0}^{\infty },$ $\left\{
I_{k}\right\} _{k=0}^{\infty }$ of $R^{n}$ and partition of unity $\left\{
\varphi _{k}\right\} _{k=0}^{\infty }$ defined e.g. in $\left[ \text{19}%
\right] .$

\textbf{Remark 2.2. }Note the following useful properties are satisfied:

supp $\varphi _{k}\subset \bar{I}_{k}$ for each $k\in \mathbb{N}_{0};$ $%
\dsum\limits_{k=0}^{\infty }\varphi _{k}\left( s\right) =1$ for each $s\in
R^{n};$ $I_{m}\cap $ supp $\varphi _{k}=\oslash $ if $\left\vert
m-k\right\vert >1;$ $\varphi _{k-1}\left( s\right) +\varphi _{k}\left(
s\right) +\varphi _{k+1}\left( s\right) =1$ for each $s\in $ supp $\varphi
_{k}$ and $k\in \mathbb{N}_{0}.$

Among the many equivalent descriptions of Besov spaces, the most useful one
for usis given in terms of the so called Littlewood-Paley decomposition.
This means that we consider $f\in S^{\prime }\left( E\right) $ as a
distributional sum $f=\dsum\limits_{k}f_{k}$ analytic functions $f_{k}$
whose Fourier transforms have support in dyadic-like $I_{k}$ and then define
the Besov norm in terms of the $f_{k}$'s.

\textbf{Definition 2.3}. Let $\gamma \in A_{q}$, $1\leq r,q\leq \infty $ and 
$s\in \mathbb{R}.$ The Besov space $B_{q,r,\gamma }^{s}\left( R^{n};E\right) 
$ is the space of all $f\in S^{\prime }\left( R^{n};E\right) $ for which 
\begin{equation}
\left\Vert f\right\Vert _{B_{q,r,\gamma }^{s}\left( R^{n};E\right)
}=\left\Vert \left\{ 2^{ks}\left( \check{\varphi}_{k}\ast f\right) \right\}
_{k=0}^{\infty }\right\Vert _{l_{r}\left( L_{q,\gamma }\left( R^{n};E\right)
\right) }=  \tag{2.2}
\end{equation}%
\begin{equation*}
\left\{ 
\begin{array}{c}
\left[ \dsum\limits_{k=0}^{\infty }2^{ksr}\left\Vert \check{\varphi}_{k}\ast
f\right\Vert _{L_{q,\gamma }\left( R^{n};E\right) }^{r}\right] ^{\frac{1}{r}}%
\text{, if }1\leq r<\infty \\ 
\sup\limits_{k\in \mathbb{N}_{0}}\left[ \dsum\limits_{k=0}^{\infty
}2^{ks}\left\Vert \check{\varphi}_{k}\ast f\right\Vert _{L_{q,\gamma }\left(
R^{n};E\right) }\right] \text{, if }r=\infty%
\end{array}%
\right.
\end{equation*}%
is finite. $B_{q,r,\gamma }^{s}\left( R^{n};E\right) $-together with the
norm in $(2.1)$, is a Banach space. $\mathring{B}_{q,r,\gamma }^{s}\left(
R^{n};E\right) $ is the closure of $S\left( R^{n};E\right) $ in $%
B_{q,r,\gamma }^{s}\left( R^{n};E\right) $ with the induced norm. In a
similar way as in $[$19, Lemma 3.2$]$ it can be shown that different choices
of $\left\{ \varphi _{k}\right\} $ lead to equivalent norms on $%
B_{q,r,\gamma }^{s}\left( R^{n};E\right) .$

Let $\Omega $ be a domain in $R^{n};$ $B_{q,r,\gamma }^{s}\left( \Omega
;E\right) $ denotes the space of restrictions to $\Omega $ of all functions
in $B_{q,r,\gamma }^{s}\left( R^{n};E\right) $ with the norm given by 
{\large 
\begin{equation*}
\left\Vert u\right\Vert _{B_{q,r,\gamma }^{s}\left( \Omega ;E\right)
}=\inf\limits_{g\in B_{q,r,\gamma }^{s}\left( R^{n};E\right) ,g\mid _{\Omega
}=u}\left\Vert g\right\Vert _{B_{q,r,\gamma }^{s}\left( R^{n};E\right) }.
\end{equation*}%
}

Let $l=\left( l_{1},l_{2},...,l_{n}\right) $, $s\in \mathbb{R}$ and $1\leq
q, $ $r\leq \infty .$ Here, $B_{q,r,\gamma }^{l,s}\left( \Omega ;E\right) $
denote a $E$-valued Sobolev-Besov weighted space of functions $u$ $\in
B_{q,\theta ,\gamma }^{s}\left( \Omega ;E\right) $ that have generalized
derivatives $D_{k}^{l_{k}}u=\frac{\partial ^{l_{k}}}{\partial x_{k}^{l_{k}}}%
u\in B_{q,r,\gamma }^{s}\left( \Omega ;E\right) ,k=1,2,...,n$ with the norm 
\begin{equation*}
\left\Vert u\right\Vert _{B_{q,\theta ,\gamma }^{l,,s}\left( \Omega
;E\right) }=\left\Vert u\right\Vert _{B_{q,r,\gamma }^{s}\left( \Omega
;E\right) }+\sum\limits_{k=1}^{n}\left\Vert D_{k}^{l_{k}}u\right\Vert
_{B_{q,\theta r\gamma }^{s}\left( \Omega ;E\right) }<\infty .
\end{equation*}

Let $E_{0}$ is continuoisly and densely belongs to $E.$ $B_{q,\theta ,\gamma
}^{l,s}\left( \Omega ;E_{0},E\right) $ denotes the space $B_{q,\theta
,\gamma }^{s}\left( \Omega ;E_{0}\right) \cap B_{q,\theta ,\gamma
}^{l,s}\left( \Omega ;E\right) $ with the norm 
\begin{equation*}
\left\Vert u\right\Vert _{B_{q,\theta ,\gamma }^{l,s}}=\left\Vert
u\right\Vert _{B_{q,\theta ,\gamma }^{l,s}\left( \Omega ;E_{0},E\right)
}=\left\Vert u\right\Vert _{B_{q,\theta ,\gamma }^{s}\left( \Omega
;E_{0}\right) }+\sum\limits_{k=1}^{n}\left\Vert D_{k}^{l_{k}}u\right\Vert
_{B_{q,\theta ,\gamma }^{s}\left( \Omega ;E\right) }<\infty .
\end{equation*}

Let $(E(X);E^{\ast }(X^{\ast }))$ be one of the pairs 
\begin{equation*}
\left( L_{q,\gamma }\left( X\right) ,L_{q^{\prime },\gamma ^{\prime }}\left(
X^{\ast }\right) \right) \text{, }\left( B_{q,r,\gamma }^{s}\left( X\right)
,B_{q^{\prime },r^{\prime },\gamma ^{\prime }}^{-s}\left( X^{\ast }\right)
\right) ,
\end{equation*}
when $1\leq q,r\leq \infty $, where%
\begin{equation*}
\gamma ^{\prime }\left( .\right) =\gamma ^{-\frac{1}{q-1}}\left( .\right) .
\end{equation*}

There is an embedding of $E^{\ast }(X^{\ast })\subset \left[ E(X)\right]
^{\ast }$\ as a norming subspace for $E(X)$. This embedding is given by the
duality map

\begin{equation*}
\langle .,.\rangle _{E(X)}:E^{\ast }(X^{\ast })\times E(X)\rightarrow 
\mathbb{C},
\end{equation*}%
where 
\begin{equation*}
\langle g,f\rangle _{L_{q,\gamma }\left( X\right)
}=\dint\limits_{R^{n}}\langle g\left( t\right) ,f\left( t\right) \rangle
_{X}dt=\dint\limits_{R^{n}}g\left( t\right) f\left( t\right) dt
\end{equation*}%
in weighted Lebesgue space setting with $E=L_{q,\gamma }$ and 
\begin{equation}
\langle g,f\rangle _{B_{q,r,\gamma }^{s}\left( X\right)
}=\dsum\limits_{n,m\in \mathbb{N}_{0}}\langle \check{\varphi}_{n}\ast g,%
\check{\varphi}_{m}\ast f\rangle _{L_{q,\gamma }\left( X\right) }  \tag{2.3}
\end{equation}%
in Besov space setting with $E=B_{q,r,\gamma }^{s}\left( X\right) .$ One can
check that this definition of duality is independent of the choice of the $%
\left\{ \varphi _{k}\right\} _{k=0}^{\infty }$.

\begin{center}
\bigskip \textbf{3.} \textbf{The Foruier transform in weighted Besov spaces}
\end{center}

\bigskip By applying the Hausdor-Young inequality we get the following
estimates for the Fourier transform on Besov spaces

\bigskip \textbf{Theorem 3.1.} Assume $\gamma \in A_{\nu }$ for $\nu \in %
\left[ 1,\infty \right] $. Let $E$ be a Banach space with weighted Fourier
type $\gamma $ and $p\in \left[ 1,2\right] .$ Let $1\leq q\leq p^{\prime }$
and $s\geq n\left( \frac{1}{q}-\frac{1}{p^{\prime }}\right) $ and $1\leq
r\leq \infty .$ Then there exists constant $C$, depending only on $%
C_{0}\left( p,\gamma ,E\right) $ so that if $f\in B_{q,r,\gamma }^{s}\left(
R^{n};E\right) $ then%
\begin{equation}
\left\Vert \left\{ \hat{f}\chi _{J_{m}}\right\} _{m=0}^{\infty }\right\Vert
_{l_{r}\left( L_{q,\gamma }\left( R^{n};E\right) \right) }\leq C\left\Vert
f\right\Vert _{B_{q,r,\gamma }^{s}\left( R^{n};E\right) },  \tag{3.1}
\end{equation}%
where $C_{0}\left( p,\gamma ,E\right) $ is a positive constant defined in
the Definition 2.1.

An immediate corollary of Theorem 3.1 follows by choosing for $q=r=1$ and $%
r=q=p^{\prime }$ we obtain respectively

\bigskip \textbf{Corollary 3.1. }Assume $\gamma \in A_{q}$ for $q\in \left[
1,\infty \right] $. Let $E$ be a Banach space with Fourier type $p\in \left[
1,2\right] .$ Then the Fourier transform $F$ defines the following bounded
operators%
\begin{equation}
F:B_{p,1,\gamma }^{\frac{n}{p}}\left( R^{n};E\right) \rightarrow L_{1,\gamma
}\left( R^{n};E\right)  \tag{3.2}
\end{equation}

\begin{equation}
F:B_{p,p^{\prime },\gamma }^{0}\left( R^{n};E\right) \rightarrow
L_{p^{\prime },\gamma }\left( R^{n};E\right) .  \tag{3.3}
\end{equation}

\bigskip The norms of the above maps $F$ are bounded above by a constant
depending only on $C_{0}\left( n,E\right) .$

\bigskip Theorem 3.1 and Corollary 3.2 remain valid if $F$ is replaced with $%
F^{-1}.$

\textbf{Proof of Theorem 3.1. }Let $f\in B_{q,r,\gamma }^{s}\left(
R^{n};E\right) .$ Then, for each $k\in \mathbb{N}_{0}$, since $\check{\varphi%
}_{k}\ast f\in L_{p,\gamma }\left( R^{n};E\right) $ and $E$ has weighted
Fourier type $\gamma $ and $p$, 
\begin{equation*}
\varphi _{k}.\hat{f}=F\left( \check{\varphi}_{k}\ast f\right) \in
L_{p^{\prime },\gamma }\left( R^{n};E\right) .
\end{equation*}

Thus by Remark 2.2, 
\begin{equation*}
\hat{f}\chi _{J_{m}}=\left( \dsum\limits_{k=m-1}^{m+1}\varphi _{k}.\hat{f}%
\right) \chi _{J_{m}}\in L_{q,\gamma }\left( R^{n};E\right) \text{ for each }%
m\in \mathbb{N}_{0}.
\end{equation*}

\bigskip Moreover, by Definition 2.2 we get 
\begin{equation*}
\left\Vert \varphi _{k}\hat{f}\right\Vert _{L_{p^{\prime },\gamma }\left(
R^{n};E\right) }=\left\Vert F\left( \check{\varphi}_{k}\ast f\right)
\right\Vert _{L_{p^{\prime },\gamma }\left( R^{n};E\right) }\leq
C_{0}\left\Vert \check{\varphi}_{k}.f\right\Vert _{L_{p,\gamma }\left(
R^{n};E\right) },
\end{equation*}%
i.e. 
\begin{equation}
\dsum\limits_{k=m-1}^{m+1}2^{ks}\left\Vert \varphi _{k}\hat{f}\right\Vert
_{L_{p^{\prime },\gamma }\left( R^{n};E\right) }\leq
C_{0}\dsum\limits_{k=m-1}^{m+1}2^{ks}\left\Vert \check{\varphi}%
_{k}.f\right\Vert _{L_{p,\gamma }\left( R^{n};E\right) }\leq  \tag{3.4}
\end{equation}

\begin{equation*}
CC_{0}\left( p,\gamma ,E\right) \left\Vert f\right\Vert _{B_{q,r,\gamma
}^{s}\left( R^{n};E\right) }.
\end{equation*}

In view of $\left( 3.4\right) $, it suffices to show that there exists the
pozitive constant $C_{1}$ so that the following holds%
\begin{equation}
\left\Vert \hat{f}\chi _{J_{m}}\right\Vert _{L_{q,\gamma }\left(
R^{n};E\right) }\leq C_{1}\dsum\limits_{k=m-1}^{m+1}2^{ks}\left\Vert \varphi
_{k}\hat{f}\right\Vert _{L_{p^{\prime },\gamma }\left( R^{n};E\right) }. 
\tag{3.5}
\end{equation}

\bigskip Firstly, consider the case where $q\neq p^{\prime }$. Choose $1\leq
\sigma <p$ that $\frac{1}{q}=\frac{1}{p^{\prime }}+\frac{1}{\sigma };$ so, $%
\frac{n}{\sigma }\leq s.$ By the generalized H\"{o}lder's inequality for
each $m\in \mathbb{N}_{0},$%
\begin{equation}
\left\Vert \hat{f}\chi _{J_{m}}\right\Vert _{L_{q,\gamma }\left(
R^{n};E\right) }\leq \dsum\limits_{k=m-1}^{m+1}\left\Vert \varphi _{k}\hat{f}%
\chi _{J_{m}}\right\Vert _{L_{q,\gamma }\left( J_{m};E\right) }\leq 
\tag{3.6}
\end{equation}%
\begin{equation*}
\dsum\limits_{k=m-1}^{m+1}\left\Vert \varphi _{k}\left( \frac{1+\left\vert
.\right\vert }{4}\right) ^{\frac{n}{\sigma }}\hat{f}\gamma ^{\frac{1}{%
p^{\prime }}}\left( .\right) \right\Vert _{L_{p^{\prime }}\left(
J_{m};E\right) }\left\Vert \gamma ^{\frac{1}{p}}\left( \frac{1+\left\vert
.\right\vert }{4}\right) ^{-\frac{n}{\sigma }}\right\Vert _{L_{\sigma
}\left( J_{m}\right) }\leq
\end{equation*}%
\begin{equation*}
\dsum\limits_{k=m-1}^{m+1}\left\Vert \hat{f}\varphi _{k}\right\Vert
_{L_{p^{\prime },\gamma }\left( J_{m};E\right) }\left\Vert \left( \frac{%
1+\left\vert .\right\vert }{4}\right) ^{\frac{n}{\sigma }}\chi
_{J_{m}}\right\Vert _{L_{\infty }}\leq
C_{2}\dsum\limits_{k=m-1}^{m+1}2^{ks}\left\Vert \hat{f}\varphi
_{k}\right\Vert _{L_{p^{\prime },\gamma }\left( J_{m};E\right) },
\end{equation*}%
where $C_{2}$ is a positive constant defined by%
\begin{equation*}
C_{2}=\left\Vert \gamma ^{\frac{1}{p}}\left( \frac{1+\left\vert .\right\vert 
}{4}\right) ^{-\frac{n}{\sigma }}\right\Vert _{L_{\sigma }\left(
J_{m}\right) }\leq \left\Vert \left( \frac{1+\left\vert .\right\vert }{4}%
\right) ^{-n}\right\Vert _{L_{\infty }\left( J_{m}\right) }\left\Vert \gamma
^{\frac{\sigma }{p}}\right\Vert _{L\left( J_{m}\right) }\leq
\end{equation*}%
\begin{equation}
4^{n}\left[ \sup_{m\in \mathbb{N}_{0}}2^{-\left( m-1\right)
n}\dint\limits_{J_{m}}\gamma ^{\frac{\sigma }{p}}\left( s\right) ds\right] ^{%
\frac{1}{\sigma }}.  \tag{3.7}
\end{equation}

Since $\gamma \in A_{\nu }$ we have 
\begin{equation*}
\sup_{m\in \mathbb{N}_{0}}2^{-\left( m-1\right) n}\dint\limits_{J_{m}}\gamma
^{\frac{\sigma }{p}}\left( s\right) ds<\infty .
\end{equation*}

For $q=p^{\prime }$ and for each $m\in \mathbb{N}$ we get 
\begin{equation*}
\left\Vert \hat{f}\chi _{J_{m}}\right\Vert _{L_{q,\gamma }\left(
R^{n};E\right) }\leq \dsum\limits_{k=m-1}^{m+1}\left\Vert \varphi _{k}\hat{f}%
\chi _{J_{m}}\right\Vert _{L_{p^{\prime },\gamma }\left( J_{m};E\right) }\leq
\end{equation*}

\begin{equation}
\dsum\limits_{k=m-1}^{m+1}2^{ks}\left\Vert \varphi _{k}\hat{f}\right\Vert
_{L_{p^{\prime },\gamma }\left( R^{n};E\right) }.  \tag{3.8}
\end{equation}%
So, from $\left( 3.6\right) $-$\left( 3.8\right) $ we obtain $\left(
3.5\right) .$

\bigskip \textbf{Remark 3.1.} By using the embedding $W_{p,\gamma
}^{j}\left( R^{n};E\right) \subset B_{q,r,\gamma }^{s}\left( R^{n};E\right) $
for $s<j\in \mathbb{N}$ we get that the statement of Theorem 3.1 remains
valid if $B_{q,r,\gamma }^{s}\left( R^{n};E\right) $ is replaced by $%
W_{p,\gamma }^{j}\left( R^{n};E\right) $.

Also, it follows from Corollary $3.2$ that if $E$ has weighted Fourier type
for $\gamma \in A_{\nu }$, $p\in \lbrack 1;2]$ and $j>\frac{n}{p}$ then the
Fourier transform $F$ defines bounded operator:%
\begin{equation*}
W_{p,\gamma }^{j}\left( R^{n};E\right) \rightarrow L_{1,\gamma }\left(
R^{n};E\right) .
\end{equation*}

Furthermore, if $E$ has weighted Fourier type for $\gamma \in A_{\nu }$, $%
p\in \lbrack 1,2]$ and $j>\frac{n}{p}$ then there is a constant $C$ so that%
\begin{equation}
\left\Vert \hat{f}\right\Vert _{L_{1,\gamma }\left( R^{n};E\right) }\leq
C\left\Vert f\right\Vert _{L_{p,\gamma }\left( R^{n};E\right) }^{1-\frac{n}{%
jp}}\left[ \dsum\limits_{\left\vert \alpha \right\vert =j}\left\Vert
D^{\alpha }f\right\Vert _{L_{p,\gamma }\left( R^{n};E\right) }\right] ^{%
\frac{n}{jp}}  \tag{3.9}
\end{equation}%
for each $f\in W_{p,\gamma }^{j}\left( R^{n};E\right) .$

\begin{center}
\bigskip \textbf{4. Fourier multipliers on weighted Lebesque spaces}
\end{center}

Consider the bounded measurable function $m:R^{n}\rightarrow B\left(
E_{1},E_{2}\right) .$ In this section, we identify conditions on $m$,
generalizing the classical Mihlin condition so that the multiplication
operator induced by $m$, i.e. the operator: $u\rightarrow T_{m}=F^{-1}mFu$
is bounded from $L_{q,\gamma }\left( R^{n};E_{1}\right) $ to $L_{q,\gamma
}\left( R^{n};E_{2}\right) .$ We will rst give rather general criteria for
Fourier multipliers in terms of the weighted Besov norm of the multiplier
function; later we derive from these results analogues of the classical
Mihlin and H\"{o}rmander conditions. To simplify the statements of our
results, we let 
\begin{equation*}
M_{p,\gamma }\left( m\right) =\inf\limits_{a>0}\left\{ \left\Vert m\left(
a,.\right) \right\Vert _{B_{p,1.\gamma }^{\frac{n}{p}}\left( R^{n};B\left(
E_{1},E_{2}\right) \right) }\right\} .
\end{equation*}

Let 
\begin{equation*}
X_{k}=L_{q,\gamma }\left( E_{k}\right) =L_{q,\gamma }\left(
R^{n};E_{k}\right) \text{, }k=1,2\text{, }Y=B_{p,1.\gamma }^{\frac{n}{p}%
}\left( R^{n};B\left( E_{1},E_{2}\right) \right) .
\end{equation*}

First we give a multiplier result from $X_{1}$ to $X_{2}$ in the spirit of
Steklin's theorem.

\textbf{Theorem 4.1. }Assume $\gamma \in A_{\nu }$ for $\nu \in \left[
1,\infty \right] $. Let $E_{1}$, $E_{2}$ be a Banach spaces with weighted
Fourier type $\gamma $ and $p\in \left[ 1,2\right] .$ Then there is a
constant $C$, depending only on $C_{01}\left( p,\gamma ,E_{1}\right) $ and $%
C_{02}\left( p,\gamma ,E_{2}\right) $, so that if $m\in Y,$ then $m$ is a
Fourier multiplier from $X_{1}$ to $X_{2}$ and 
\begin{equation*}
\left\Vert T_{m}\right\Vert _{B\left( X_{1},X_{2}\right) }\leq CM_{p,\gamma
}\left( m\right)
\end{equation*}%
for each $q\in \left[ 1,\infty \right] $.

\bigskip Let $E^{\ast }$ denotes the dual space of $E$ and $A^{\ast }-$%
denotes the conjugate of the operator $A.$

The proof of Theorem 4.3 uses the following lemma.

\bigskip \textbf{Lemma 4.1. }Assume $\gamma \in A_{q}$ for $q\in \left[
1,\infty \right] $ and $k\in L_{1}\left( R^{n}\text{; }B\left(
E_{1},E_{2}\right) \right) .$ Suppose that there exists constants $C_{i}$ so
that for each $x\in E_{1}$ and $x^{\ast }\in E_{2}^{\ast }$%
\begin{equation}
\dint\limits_{R^{n}}\left\Vert k\left( s\right) x\right\Vert _{E_{2}}ds\leq
M_{0}\left\Vert x\right\Vert _{E_{1}},\text{ }\dint\limits_{R^{n}}\left\Vert
k^{\ast }\left( s\right) x^{\ast }\right\Vert _{E_{1}^{\ast }}ds\leq
M_{1}\left\Vert x^{\ast }\right\Vert _{E_{2}^{\ast }}.  \tag{4.1}
\end{equation}

\bigskip Then the convolution operator $K:$ $X_{1}\rightarrow X_{2}$ defined
by 
\begin{equation}
\left( Kf\right) \left( t\right) =\dint\limits_{R^{n}}k\left( t-s\right)
f\left( s\right) ds\text{ for }t\in R^{n}  \tag{4.2}
\end{equation}%
satisfies that 
\begin{equation*}
\left\Vert K\right\Vert _{B\left( X_{1},X_{2}\right) }\leq M_{0}^{\frac{1}{q}%
}M_{1}^{1-\frac{1}{q}}.
\end{equation*}

\textbf{Proof. }\ Since $k\in L_{1}\left( R^{n}\text{; }B\left(
E_{1},E_{2}\right) \right) $ it is well-known that $\left( 4.2\right) $
defines a bounded operator on $X_{1}.$ Indeed, for $f\in X_{1}\cap L_{\infty
}\left( R^{n};E_{1}\right) $ we have 
\begin{equation}
\dint\limits_{R^{n}}\left\Vert k\left( t-s\right) f\left( s\right)
\right\Vert _{E_{2}}ds=\dint\limits_{R^{n}}\left\Vert k\left( s\right)
f_{s}\left( t\right) \right\Vert _{E_{2}}ds\leq \left\Vert k\right\Vert
_{L_{1}\left( R^{n};B\left( E_{1}E_{2}\right) \right) }\left\Vert
f\right\Vert _{L_{\infty }\left( R^{n};E_{1}\right) }  \tag{4.3}
\end{equation}%
for each $t\in R^{n}$ and $f_{s}\left( t\right) =f\left( t-s\right) .$ From $%
\left( 4.3\right) $ by applying the Minkowski's inequality for integral with
weight $\left[ \text{20, \S\ A.1}\right] $ we get 
\begin{equation*}
\left\Vert Kf\left( .\right) \right\Vert _{X_{2}}\leq
\dint\limits_{R^{n}}\left\Vert k\left( s\right) f_{s}\left( t\right)
\right\Vert _{X_{2}}ds\leq \dint\limits_{R^{n}}\left\Vert k\left( s\right)
\right\Vert _{B\left( E_{1},E_{2}\right) }\left\Vert f_{s}\right\Vert
_{X_{1}}ds=
\end{equation*}%
\begin{equation*}
\left\Vert k\right\Vert _{L_{1}\left( R^{n}\text{; }B\left(
E_{1},E_{2}\right) \right) }\left\Vert f_{s}\right\Vert _{X_{1}}.
\end{equation*}

Now, for $q=1$ we have from $\left( 4.1\right) $ 
\begin{equation*}
\left\Vert Kf\right\Vert _{L_{1,\gamma }\left( R^{n};E_{1}\right) }\leq
\dint\limits_{R^{n}}\left( \dint\limits_{R^{n}}\left\Vert k\left( s\right)
f_{s}\left( t\right) \right\Vert _{E_{1}}ds\right) \gamma \left( t\right)
dt\leq
\end{equation*}%
\begin{equation*}
M_{0}\dint\limits_{R^{n}}\left\Vert f\left( t\right) \right\Vert
_{E_{1}}\gamma \left( t\right) dt=M_{0}\left\Vert f\right\Vert _{L_{1,\gamma
}\left( R^{n};E_{1}\right) }.
\end{equation*}

\bigskip Hence, 
\begin{equation}
\left\Vert K\right\Vert _{B\left( L_{1,\gamma }\left( R^{n};E_{1}\right)
\right) }\leq M_{0}.  \tag{4.4}
\end{equation}

If $q=\infty $, then for each $L_{\infty ,\gamma }\left( R^{n};E\right) $, $%
x^{\ast }\in E_{2}^{\ast }$ and $t\in R^{n}$ by using $\left( 4.1\right) $
we get%
\begin{equation*}
\left\vert \langle x^{\ast },\left( Kf\right) \left( t\right) \rangle
_{E_{2}}\right\vert \leq \dint\limits_{R^{n}}\left\vert \langle k^{\ast
}\left( t-s\right) x^{\ast },f\left( s\right) \rangle _{E_{1}}\right\vert
\gamma \left( s\right) ds\leq
\end{equation*}%
\begin{equation*}
\dint\limits_{R^{n}}\left\Vert k^{\ast }\left( t-s\right) x^{\ast
}\right\Vert _{E_{1}^{\ast }}\left\Vert f\left( s\right) \right\Vert
_{E_{1}}\gamma \left( s\right) ds\leq M_{1}\left\Vert x^{\ast }\right\Vert
_{E_{1}^{\ast }}\left\Vert f\right\Vert _{L_{\infty ,\gamma }\left(
R^{n};E\right) }.
\end{equation*}

Thus, 
\begin{equation}
\left\Vert K\right\Vert _{B\left( L_{\infty ,\gamma }\left(
R^{n};E_{1}\right) \right) }\leq M_{1}.  \tag{4.5}
\end{equation}

Let $L_{\infty ,\gamma }\left( R^{n};E_{1}\right) $ denotes the closure in $%
L_{\infty ,\gamma }\left( R^{n};E_{1}\right) $ norm of the simple functions $%
\dsum\limits_{k=1}^{m}x_{k}\chi _{A_{k}},$ where $x_{k}\in E_{1}$, vol $%
A_{k}<\infty $ and $m\in \mathbb{N}.$ Then one can check that $K$ maps $%
L_{\infty ,\gamma }\left( R^{n};E_{1}\right) $ into $L_{\infty ,\gamma
}\left( R^{n};E_{2}\right) $. Indeed, for $f=\chi _{A}$, we have 
\begin{equation*}
Kf\left( t\right) =\dint\limits_{t-A}k\left( s\right) xds\rightarrow 0\text{
for }t\rightarrow \infty
\end{equation*}%
and $Kf$ is a continuous function from $R^{n}$ to $E_{2}.$ Now, the
Riesz-Thorin theorem (cf. [5, Thm 5.1.2]) yields the claim for $1<p<\infty .$

\textbf{Proof of Theorem 4.1. }\ First assume in addition that $m\in S\left(
B\left( E_{1},E_{2}\right) \right) .$ Hence, $\check{m}\in S\left( B\left(
E_{1},E_{2}\right) \right) .$ Fix $x\in E_{1}$. For an appropriate choice of 
$a>0$, we can apply Corollary 3.1 to the function $t\rightarrow m\left(
at\right) x$ in $B_{p,1.\gamma }^{\frac{n}{p}}\left( R^{n};E_{2}\right) $
and use that%
\begin{equation*}
F^{-1}\left[ m\left( a.\right) x\right] \left( s\right) =a^{-n}\check{m}%
\left( \frac{s}{a}\right) x
\end{equation*}%
to get 
\begin{equation*}
\left\Vert \check{m}\left( .\right) x\right\Vert _{L_{1,\gamma }\left(
R^{n};E_{1}\right) }=\left\Vert F^{-1}m\left( a.\right) x\right\Vert
_{L_{1,\gamma }\left( R^{n};E_{1}\right) }\leq
\end{equation*}%
\begin{equation*}
C_{1}\left\Vert m\left( a.\right) x\right\Vert _{B_{p,1,\gamma }^{\frac{n}{p}%
}\left( R^{n};B\left( E_{1},E_{2}\right) \right) }\left\Vert x\right\Vert
_{E_{1}}\leq 2C_{1}M_{p,\gamma }\left\Vert x\right\Vert _{E_{1}},
\end{equation*}%
for some constant $C_{1}$ which depends on $C_{0}\left( p,\gamma
,E_{2}\right) .$

By the additional assumption on $m$ we get 
\begin{equation*}
m^{\ast }\left( .\right) \in S\left( B\left( E_{2}^{\ast },E_{1}^{\ast
}\right) \right) \text{, and }F^{-1}m^{\ast }\left( .\right) =\left[ \check{m%
}\left( .\right) \right] ^{\ast }\in S\left( B\left( E_{2}^{\ast
},E_{1}^{\ast }\right) \right) .
\end{equation*}

Let $x^{\ast }\in E_{2}^{\ast }.$ Similarly, by applying Corollary 3.1 to an
appropriate function%
\begin{equation*}
t\rightarrow \left[ m\left( at\right) \right] ^{\ast }x^{\ast }\text{ in }%
B_{p,1.\gamma }^{\frac{n}{p}}\left( R^{n};E_{1}^{\ast }\right)
\end{equation*}%
and using the fact that \ $M_{p,\gamma }\left( m\right) =M_{p,\gamma }\left(
m^{\ast }\right) $, one has%
\begin{equation*}
\left\Vert \left[ \check{m}\left( .\right) \right] ^{\ast }x^{\ast
}\right\Vert _{L_{1,\gamma }\left( R^{n};E_{1}^{\ast }\right) }\leq
2C_{2}M_{p,\gamma }\left( m\right) \left\Vert x^{\ast }\right\Vert
_{E_{2}^{\ast }}
\end{equation*}%
for some constant $C_{2}$ which depends $C_{0}\left( p,\gamma ,E_{1}^{\ast
}\right) .$ By Lemma 4.1, the convolution operator%
\begin{equation*}
\left( T_{m}f\right) \left( t\right) =\dint\limits_{R^{n}}\check{m}\left(
t-s\right) f\left( s\right) ds
\end{equation*}%
satisfies 
\begin{equation*}
\left\Vert T_{m}\right\Vert _{B\left( X_{1},X_{2}\right) }\leq CM_{p,\gamma
}\left( m\right) ,
\end{equation*}%
where $C=2\max \left\{ C_{1},C_{2}\right\} .$ Furthermore, since $m\in
L_{1}\left( R^{n};B\left( E_{1},E_{2}\right) \right) $, then $T_{m}$
satisfies the following%
\begin{equation}
T_{m}f=F^{-1}m\left( .\right) f\left( .\right) \text{ for all }f\in S\left(
R^{n};E_{1}\right) ,  \tag{4.6}
\end{equation}%
also%
\begin{equation}
T_{m}\in C\left( \sigma \left( X_{1},X_{1}^{\ast }\right) ,\sigma \left(
X_{2},X_{2}^{\ast }\right) \right) ,  \tag{4.7}
\end{equation}%
where $\sigma \left( X_{k},X_{k}^{\ast }\right) $ denote the interpolation
spaces of $X_{k},$ $X_{k}^{\ast }.$

For the general case, let $m\in Y.$ It is known that $S\left( R^{n};B\left(
E_{1},E_{2}\right) \right) $ is dence in $Y$ when $\gamma \in A_{\nu },$ $%
\nu \in \left[ 1,\infty \right] .$ Now, let we choose a sequence $\left\{
m_{n}\right\} _{n}^{\infty }\subset S\left( R^{n};B\left( E_{1},E_{2}\right)
\right) $ that converges to $m$ in the $Y-$norm and obtain operators $%
T_{m_{n}}\in B\left( X_{1},X_{2}\right) $, where 
\begin{equation*}
T_{m_{n}}f=F^{-1}m_{n}\left( .\right) f\left( .\right) ,\text{ }f\in X_{1}.
\end{equation*}

It is clear to see that, the properties $\left( 4.6\right) $ and $\left(
4.7\right) $ pass from $T_{m_{n}}$ to $T_{m}.$ One also has that%
\begin{equation*}
\left\Vert T_{m}\right\Vert _{B\left( X_{1},X_{2}\right) }\leq C\left\Vert
m\right\Vert _{Y}\text{.}
\end{equation*}

Fix $a>0$ such that $m\left( a.\right) \in Y.$ Then $I_{E_{2}}\circ
T_{m\left( a.\right) }=T_{m}\circ I_{E_{1}}$, where $I_{\mathbb{Z}%
}:L_{q,\gamma }\left( R^{n};\mathbb{Z}\right) \rightarrow L_{q,\gamma
}\left( R^{n};\mathbb{Z}\right) $ is the isometry 
\begin{equation*}
T\left( f\right) \left( t\right) =a^{\frac{n}{q}}f\left( at\right) .
\end{equation*}%
Thus, 
\begin{equation*}
\left\Vert T_{m}\right\Vert _{B\left( X_{1},X_{2}\right) }=\left\Vert
T_{m\left( a.\right) }\right\Vert _{B\left( X_{1},X_{2}\right) }\leq
C\left\Vert m\right\Vert _{Y},
\end{equation*}%
i.e. 
\begin{equation*}
\left\Vert T_{m}\right\Vert _{B\left( X_{1},X_{2}\right) }\leq CM_{p,\gamma
}\left( m\right) .
\end{equation*}

\bigskip The following remark collects some basic facts about the Fourier
multiplier operators $T_{m}$ given in Theorem 4.1 that will be used in the
proof of Theorem 4.2.

\textbf{Remark 4.1. }Let $f$ $\in X_{1}$ and $\Omega $ be a closed subset of 
$R^{n}.$ Then the following are valid:

(a) Viewing $f$ and $T_{m}f$ as distributions, if supp $\hat{f}\subset
\Omega $ then supp $F\left( T_{m}f\right) \subset \Omega ;$

(b) $T_{m_{1}+m_{2}}=T_{m_{1}}+T_{m_{2}}.$ If $\varphi \in S$, then $\check{%
\varphi}\ast T_{m}f=T_{m}\left( \check{\varphi}\ast f\right) =T_{\varphi
m}\left( f\right) ;$

(c) If $\varphi \in S$ is $1$ on supp $\hat{f}$, then $T_{\varphi m}\left(
f\right) =T_{m}\left( f\right) ;$

(d) $T_{m}^{\ast }$\ \ restricted to $L_{q`,\gamma }\left( R^{n};E_{2}^{\ast
}\right) $ is $T_{m^{\ast }\left( -.\right) }.$

\begin{center}
\bigskip \textbf{5. Fourier multipliers on weighted Besov spaces}
\end{center}

Consider the bounded measurable function $m:R^{n}\rightarrow B\left(
E_{1},E_{2}\right) .$ In this section we identify conditions on $m$,
generalizing the classical Mikhlin condition so that the multiplication
operator induced by $m$, i.e. the operator: $u\rightarrow T_{m}=F^{-1}mFu$
is bounded from $B_{p,q,\gamma }^{s}\left( R^{n};E_{1}\right) $ to $%
B_{p,q,\gamma }^{s}\left( R^{n};E_{2}\right) .$

\bigskip By applying this Theorem 4.1 to the blocks of the Littlewood Paley
decomposition of Besov spaces we will now get the main result of this
section. Let 
\begin{equation*}
Y_{i}=B_{q,r,\gamma }^{s}\left( R^{n};E_{i}\right) \text{, }i=1,2.
\end{equation*}

\bigskip \textbf{Theorem 5.1. }Assume $\gamma \in A_{\nu }$ for $\nu \in %
\left[ 1,\infty \right] $. Let $E_{1}$, $E_{2}$ be a Banach spaces with
weighted Fourier type $\gamma $ and $p\in \left[ 1,2\right] .$ Then there is
a constant $C$ depending only on $C_{01}\left( p,\gamma ,E_{1}\right) $ and $%
C_{02}\left( p,\gamma ,E_{2}\right) $, so that if 
\begin{equation}
\varphi _{k}m\in Y\text{ and }M_{p,\gamma }\left( \varphi _{k}m\right) \leq A%
\text{ for each }k\in \mathbb{N}_{0}  \tag{5.1}
\end{equation}%
then $m$ is a Fourier multiplier from $Y_{1}$ to $Y_{2}$ and 
\begin{equation*}
\left\Vert T_{m}\right\Vert _{B\left( Y_{1},Y_{2}\right) }\leq CA
\end{equation*}%
for each $s\in \mathbb{R}$ and $q,$ $r\in \left[ 1,\infty \right] .$

\textbf{Proof. }By definition partition of unity $\left\{ \varphi
_{k}\right\} _{k=0}^{\infty }$ we have 
\begin{equation*}
T_{m}f=F^{-1}m\hat{f}=\dsum\limits_{k\in \mathbb{N}_{0}}F^{-1}\left[ \left(
\varphi _{k-1}+\varphi _{k}+\varphi _{k+1}\right) mF\left[ \left( \check{%
\varphi}_{k}\ast f\right) \right] \right] =
\end{equation*}%
\begin{equation}
\dsum\limits_{k\in \mathbb{N}_{0}}T_{\left( \varphi _{k-1}+\varphi
_{k}+\varphi _{k+1}\right) m}\left( \check{\varphi}_{k}\ast f\right) \text{,}
\tag{5.2}
\end{equation}%
where $T_{m}$ is the Fourier multiplier operator on $X_{1}$\ given by
Theorem 4.1. Theorem 4.1 gives that $m\varphi _{k}$ induces a Fourier
multiplier operator $T_{m.\varphi _{k}}$\ with%
\begin{equation*}
\left\Vert T_{m.\varphi _{k}}\right\Vert _{B\left( X_{1},X_{2}\right) }\leq
CM_{p,\gamma }\left( \varphi _{k}m\right) \leq CA
\end{equation*}%
for some constant $C$ depending only on $C_{0,1}\left( p,\gamma
,E_{1}\right) $ and $C_{0,2}\left( p,\gamma ,E_{2}\right) .$ Let%
\begin{equation*}
\psi _{k}=\varphi _{k-1}+\varphi _{k}+\varphi _{k+1}.
\end{equation*}

Note that $\psi _{k}\left( s\right) \equiv 1$ when $s\in $ supp $\varphi
_{k} $. Then $m\psi _{k}$ induces the Fourier multiplier operator $T_{m.\psi
_{k}} $ with 
\begin{equation*}
T_{m\psi _{k}}=T_{m\varphi _{k-1}}+T_{m\varphi _{k}}+T_{m\varphi _{k+1}}\in
B\left( X_{1},X_{2}\right)
\end{equation*}%
and 
\begin{equation*}
\left\Vert T_{m.\psi _{k}}\right\Vert _{B\left( X_{1},X_{2}\right) }\leq 3CA.
\end{equation*}%
Define $T_{0}$: $S(E_{1})\rightarrow S^{\prime }(E_{1})$ by 
\begin{equation*}
T_{0}f=F^{-1}m\left( .\right) Ff\left( .\right) .
\end{equation*}%
If $f\in S(E_{1}),$ then%
\begin{equation*}
\check{\varphi}_{k}\ast T_{0}f=T_{m\psi _{k}}\left( \check{\varphi}_{k}\ast
f\right)
\end{equation*}%
for each $k\in \mathbb{N}_{0}$ since 
\begin{equation*}
F\left[ T_{m\psi _{k}}\left( \check{\varphi}_{k}\ast f\right) \right] \left(
.\right) =m\left( .\right) \psi _{k}\left( .\right) F\left[ \left( \check{%
\varphi}_{k}\ast f\right) \left( .\right) \right] =
\end{equation*}%
\begin{equation*}
\varphi _{k}\left( .\right) m\left( .\right) \hat{f}\left( .\right) =\varphi
_{k}\left( .\right) F\left( T_{0}f\right) =F\left[ \left( \check{\varphi}%
_{k}\ast T_{0}f\right) \left( .\right) \right] .
\end{equation*}

So, by the definition of the Besov norm 
\begin{equation*}
\left\Vert T_{0}f\right\Vert _{Y_{2}}\leq 3CA\left\Vert T_{0}f\right\Vert
_{Y_{1}}.
\end{equation*}

Thus $T_{0}$ extends to a bounded linear operator from $\mathring{B}%
_{q,r,\gamma }^{s}\left( R^{n};E_{1}\right) $ to 
\begin{equation*}
\mathring{B}_{q,r,\gamma }^{s}\left( R^{n};E_{2}\right) .
\end{equation*}
If $q,$ $r<\infty $ then $\mathring{B}_{q,r,\gamma }^{s}\left(
R^{n};E\right) =$ $B_{q,r,\gamma }^{s}\left( R^{n};E\right) $ and so all
that would remain is to verify the weak continuity condition $(4.7$).
However, we continue with the proof in order to also cover the case $%
q=\infty $ or$\ r=\infty .$ We shall show that the operator $%
T_{m}:Y_{1}\rightarrow Y_{2}$ defined by 
\begin{equation}
T_{m}f=\dsum\limits_{k=1}^{\infty }f_{k}\text{, }f_{k}=T_{m\psi _{k}}\left( 
\check{\varphi}_{k}\ast f\right) \in X_{2}  \tag{5.3}
\end{equation}%
is indeed a (norm) continuous operator. Fix $f\in Y_{1}.$ First we show that
the formal series $(5.3)$ defines an element in $S^{\prime }(E_{2}).$
Towards this, fix $\varphi \in S$. Remark $4.1$ gives that supp $%
f_{k}\subset \bar{I}_{k}$. Thus%
\begin{equation*}
f_{k}\left( \varphi \right) =\hat{f}_{k}\left( \check{\varphi}\right) =\hat{f%
}_{k}\left( \psi _{k}\left( -.\right) \check{\varphi}\right) =f_{k}\left(
\psi _{k}\ast \varphi \right)
\end{equation*}%
and so by using H\"{o}lder inequality with weight $\gamma \in A_{q}$ as in $%
\left( 3.7\right) $ we get%
\begin{equation*}
\dsum\limits_{k=1}^{\infty }\left\Vert f_{k}\left( \varphi \right)
\right\Vert _{E_{2}}\leq \dsum\limits_{k=1}^{\infty }\left\Vert
f_{k}\right\Vert _{X_{2}}\left\Vert \gamma ^{-\frac{1}{q}}\left( \psi
_{k}\ast \varphi \right) \right\Vert _{L_{q^{\prime }}\left( \mathbb{C}%
\right) }\leq
\end{equation*}%
\begin{equation*}
M\dsum\limits_{k=1}^{\infty }2^{ks}\left\Vert \check{\varphi}_{k}\ast
f\right\Vert _{X_{2}}\left\Vert 2^{-ks}\psi _{k}\ast \varphi \right\Vert
_{L_{q^{\prime },\sigma }\left( \mathbb{C}\right) }\leq
\end{equation*}%
\begin{equation*}
M2^{\left\vert s\right\vert }\left\Vert f\right\Vert _{Y_{2}}\left\Vert
\varphi \right\Vert _{B_{q^{\prime },r^{\prime },\sigma }^{-s}}\left( 
\mathbb{C}\right) ,\text{ }
\end{equation*}%
where 
\begin{equation*}
\sigma \left( .\right) =\gamma ^{1-q}\left( .\right) .
\end{equation*}%
Thus $\left( T_{m}f\right) \left( \varphi \right) $ for $\varphi \in S$
defines a linear map from $S$ into $E_{2}$ which is continuous by well known
inlusion%
\begin{equation*}
S\left( E_{2}\right) \subset Y_{2}\subset S^{\prime }\left( E_{2}\right) .
\end{equation*}%
By Remark 4.1, for each $j$, $k\in \mathbb{N}_{0}$%
\begin{equation*}
\check{\varphi}_{j}\ast T_{m\psi _{k}}\left( \check{\varphi}_{k}\ast
f\right) =T_{m\psi _{k}}\left( \check{\varphi}_{j}\ast \check{\varphi}%
_{k}\ast f\right) =\check{\varphi}_{k}\ast T_{m\psi _{k}}\left( \check{%
\varphi}_{j}\ast f\right) .
\end{equation*}%
Thus, since the support of $\varphi _{k}$ intersects the support of $\varphi
_{j}$ only for $\left\vert k-j\right\vert \leq 1$, applyin Remark 4.1
further gives%
\begin{equation}
\check{\varphi}_{k}\ast T_{m}f=\dsum\limits_{j=k-1}^{k+1}\check{\varphi}%
_{k}\ast T_{m\psi _{j}}\left( \check{\varphi}_{j}\ast f\right)
=\dsum\limits_{j=k-1}^{k+1}\check{\varphi}_{j}\ast T_{m\psi _{j}}\left( 
\check{\varphi}_{k}\ast f\right) =  \tag{5.4}
\end{equation}%
\begin{equation*}
\dsum\limits_{j=k-1}^{k+1}T_{m\varphi _{j}\psi _{j}}\left( \check{\varphi}%
_{k}\ast f\right) =T_{m\psi _{k}}\left( \check{\varphi}_{k}\ast f\right) .
\end{equation*}%
Hence, $\check{\varphi}_{k}\ast T_{m}f\in X_{2}$ and 
\begin{equation*}
\left\Vert \check{\varphi}_{k}\ast T_{m}f\right\Vert _{X_{2}}\leq
3CA\left\Vert \check{\varphi}_{k}\ast f\right\Vert _{X_{1}},
\end{equation*}%
from which and in view of $\left( 5.2\right) $ it follows that range of $%
T_{m}$ is contained in $Y_{1}$ and that norm of $T_{m}$ as an operator from $%
Y_{1}$ to $Y_{2}$ is bounded by a constant depending on the items claimed.
Furthermore, $T_{m}$ extends $T_{0}$; indeed, if $f\in S(E_{1})$ then 
\begin{equation*}
F\left( T_{m}f\right) =\dsum\limits_{k=1}^{\infty }F\left[ T_{m\psi
_{k}}\left( \check{\varphi}_{k}\ast f\right) \right] =\dsum\limits_{k=1}^{%
\infty }m\psi _{k}\varphi _{k}\hat{f}=
\end{equation*}%
\begin{equation*}
\dsum\limits_{k=1}^{\infty }m\varphi _{k}\hat{f}=F\left( T_{0}f\right) .
\end{equation*}%
It remains to show only that $T_{m}$ satisfies $(4.7)$. Since $[m(-.)]^{\ast
}$ : $R^{n}\rightarrow $ $B(E_{2}^{\ast };E_{1}^{\ast })$ also satisfies
condition $\left( 5.1\right) $, the Fourier multiplier operator $T_{m^{\ast
}(-.)}$, defined by $(4.6)$, extends to $T_{m^{\ast }(-.)}\in B\left(
E^{\ast }\left( E_{2}^{\ast }\right) ,E^{\ast }\left( E_{1}^{\ast }\right)
\right) $, for $E=B_{q,r,\gamma }^{s}.$ It suffices to show that $%
T_{m}^{\ast }$ restricted to $E^{\ast }\left( E_{2}^{\ast }\right) $ is $%
T_{m^{\ast }(-.)}.$ Hence, fix $g\in E^{\ast }\left( E_{2}^{\ast }\right) $, 
$f\in B_{q,r,\gamma }^{s}\left( E_{1}\right) $ and by using $\left(
5.4\right) $ and $\left( 2.3\right) $ we have%
\begin{equation*}
\langle T_{m}^{\ast }g,f\rangle _{Y_{1}}=\dsum\limits_{n,k\in \mathbb{N}%
_{0}}\langle \check{\varphi}_{n}\ast g,\check{\varphi}_{k}\ast T_{m}f\rangle
_{L_{q,\gamma }\left( E_{2}\right) }=
\end{equation*}%
\begin{equation}
\dsum\limits_{n,k\in \mathbb{N}_{0}}\langle \check{\varphi}_{n}\ast
g,T_{m\psi _{k}}\left( \check{\varphi}_{k}\ast f\right) \rangle
_{L_{q,\gamma }\left( E_{2}\right) }.  \tag{5.5}
\end{equation}%
and 
\begin{equation*}
\langle T_{m^{\ast }(-.)}g,f\rangle _{Y_{1}}=\dsum\limits_{n,k\in \mathbb{N}%
_{0}}\langle \check{\varphi}_{n}\ast T_{m^{\ast }(-.)}g,\check{\varphi}%
_{k}\ast f\rangle _{L_{q,\gamma }\left( E_{1}\right) }=
\end{equation*}%
\begin{equation}
\dsum\limits_{n,k\in \mathbb{N}_{0}}\langle T_{m^{\ast }(-.)\psi _{n}\left(
.\right) }\left( \check{\varphi}_{n}\ast g\right) ,\check{\varphi}_{k}\ast
f\rangle _{L_{q,\gamma }\left( E_{1}\right) }.  \tag{5.6}
\end{equation}

Fix $K_{0}\subset $ $\mathbb{N}_{0}$ and choose a radial $\psi \in S$ with
compact support such that $\psi $ is 1 on $\cup _{k=1}^{K_{0}+1}$supp $%
\varphi _{k}$. If $n$, $k\in \left\{ 0,1,...,K_{0}\right\} $, then by Remark
4.1 we get 
\begin{equation}
T_{m\psi _{k}}\left( \check{\varphi}_{k}\ast f\right) =T_{m\psi \psi
_{k}}\left( \check{\varphi}_{k}\ast f\right) =T_{m\psi }\left( \check{\varphi%
}_{k}\ast f\right)  \tag{5.7}
\end{equation}%
and 
\begin{equation}
T_{m^{\ast }(-.)\psi _{n}\left( .\right) }\left( \check{\varphi}_{n}\ast
f\right) =T_{m^{\ast }(-.)\psi \left( .\right) \psi _{n}\left( .\right)
}\left( \check{\varphi}_{n}\ast f\right) =T_{m^{\ast }(-.)\psi _{n}\left(
.\right) }\left( \check{\varphi}_{n}\ast f\right) .  \tag{5.8}
\end{equation}%
since $m\psi $ and $m^{\ast }(-.)\psi _{n}\left( .\right) $ satisfy the
assumptions of Theorem 4.1. Hence, by $\left( 5.5\right) -\left( 5.8\right) $
and by Remark 4.1 we have 
\begin{equation*}
\langle T_{m}^{\ast }g,f\rangle =\langle T_{m^{\ast }(-.)}g,f\rangle .
\end{equation*}

The next lemma gives a convenient way to verify the assumption of Theorem
4.8 in terms of derivatives.

By reasoning as Lemma 4.10 and Corollary 4.11 in $\left[ \text{11}\right] $
we obtain

\textbf{Lemma 5.1. }\ Let $\frac{n}{p}<l\in \mathbb{N}$ and $\sigma \in %
\left[ p,\infty \right] $. If $m\in C^{l}\left( R^{n};B\left(
E_{1},E_{2}\right) \right) $ and there exists a positive constant $A$ so
that 
\begin{equation}
\left\Vert D^{\alpha }m\right\Vert _{L_{\sigma }\left( R^{n};B\left(
E_{1},E_{2}\right) \right) }\leq A  \tag{5.9}
\end{equation}%
for each $k\in \mathbb{N}$, $\alpha \in \mathbb{N}_{0}^{n}$ with $\left\vert
\alpha \right\vert \leq l-1.$ Then $m$ satisfies condition $(5.1)$ of
Theorem 5.1.

\textbf{Corollary} \textbf{5.1}. Let $q,$ $r$ $\in \left[ 1,\infty \right] $
and $s\in R$. If $m\in C^{l}\left( R^{n};B\left( E_{1},E_{2}\right) \right) $
and there exists a positive constant $A$ so that 
\begin{equation}
\sup\limits_{t\in R^{n}}\left( 1+\left\vert t\right\vert \right)
^{\left\vert \alpha \right\vert }\left\Vert D^{\alpha }m\right\Vert
_{L_{\sigma }\left( R^{n};B\left( E_{1},E_{2}\right) \right) }\leq A 
\tag{5.10}
\end{equation}%
for each $k\in \mathbb{N}$, $\alpha \in \mathbb{N}_{0}^{n}$ with $\left\vert
\alpha \right\vert \leq l$ and $m_{k}\left( .\right) =m\left(
2^{k-1}.\right) .$ Then $m$\ is a Fourier multiplier from $Y_{1}$ to $Y_{2}$
provided one of the following conditions hold:

(a) $E_{1}$ and $E_{2}$ are arbitrary Banach spaces and $l=n+1;$

(b) $E_{1}$ and $E_{2}$ are uniformly convex Banach spaces and $l=n;$

(c) $E_{1}$ and $E_{2}$ have Fourier type $p$ and $l=\left[ \frac{n}{p}%
\right] +1.$

\begin{center}
\bigskip \textbf{6. Embeding theorems in Besov-Lions type spaces}
\end{center}

From $\left[ 23\right] $ we have

\textbf{Lemma 6.1}. Let $A$ be a positive operator on a Banach space\ $E$, $%
b $ be a nonnegative real number and \ $r=\left(
r_{1},r_{2},...,r_{n}\right) $ where $r_{k}\in \left\{ 0,b\right\} .$ Let $%
t=\left( t_{1},t_{2},...,t_{n}\right) $, $0<$ $t_{k}\leq T<\infty ,$ $%
k=1,2,...,n$, $\alpha =\left( \alpha _{1},\alpha _{2},...,\alpha _{n}\right) 
$ and $l=\left( l_{1},l_{2},...,l_{n}\right) $, where $l_{k}$ are positive\
and $\alpha _{k\text{ }}$are nonnegative integers such that $\varkappa
=\left\vert \left( \alpha +r\right) :l\right\vert \leq 1.$ For $0<$\ $h\leq
h_{0}<\infty $ and, $0\leq \mu \leq 1-\varkappa $ \ the operator-function \ 

\begin{equation*}
\Psi _{t}\left( \xi \right) =\Psi _{t,h,\mu }\left( \xi \right)
=\prod\limits_{k=1}^{n}t_{k}^{\frac{\alpha _{k}+r_{k}}{l_{k}}}\xi ^{r}\left(
i\xi \right) ^{\alpha }A^{1-\varkappa -\mu }h^{-\mu }\left[ A+\eta \left(
t,\xi \right) \right] ^{-1}
\end{equation*}%
is bounded operator in $E$ uniformly with respect to $\xi \in R^{n},$ $h>0$
and $t,$ i.e there is a constant $C_{\mu }$ such that 
\begin{equation*}
\ \ \ \left\Vert \Psi _{t,h,\mu }\left( \xi \right) \right\Vert _{L\left(
E\right) }\leq C_{\mu }
\end{equation*}%
for all $\xi \in R^{n}$ and $h>0,$ where, 
\begin{equation*}
\eta \left( t,\xi \right) =\sum\limits_{k=1}^{n}t_{k}\left\vert \xi
_{k}\right\vert ^{l_{k}}+h^{-1}.
\end{equation*}

Let 
\begin{equation*}
\alpha =\left( \alpha _{1},\alpha _{2},...,\alpha _{n}\right) ,\text{ }%
l=\left( l_{1},l_{2},...,l_{n}\right) \text{, }\varkappa
=\sum\limits_{k=1}^{n}\frac{\alpha _{k}}{l_{k}}.
\end{equation*}

Let $l=\left( l_{1},l_{2},...,l_{n}\right) $, where $l_{k}$ are positive
integers. Let%
\begin{equation*}
\nu \left( l\right) =\max\limits_{k,,j\in \left\{ 1,2,...,n\right\} }\left[ 
\frac{1}{l_{k}}-\frac{1}{l_{j}}\right] \text{, }\eta \left( t\right)
=\prod\limits_{k=1}^{n}t_{k}^{\frac{\alpha _{k}}{l_{k}}},\text{ }%
Y=B_{p,\theta ,\gamma }^{l,s}\left( R^{n};E\left( A\right) ,E\right) .
\end{equation*}

\ \textbf{Theorem 6.1}. Suppose the following conditions hold:

(1) $\gamma \in A_{\nu }$ for $\nu \in \left[ 1,\infty \right] $. $E$ is a
Banach spaces with weighted Fourier type $\gamma $ and $\sigma \in \left[ 1,2%
\right] $;

(2) $t=\left( t_{1},t_{2},...,t_{n}\right) ,$ $0<$ $t_{k}\leq T<\infty ,$ $%
k=1,2,...,n$, $1<p\leq q<\infty $, $\theta \in \left[ 1,\infty \right] $; $\
\ \ \ \ \ $

$\ $(3) $l_{k}$ are positive\ and $\alpha _{k\text{ }}$are nonnegative
integers such that $0<\varkappa +\nu \left( l\right) \leq 1,$ and let \ $%
0\leq \mu \leq 1-\varkappa -\nu \left( l\right) $;

(4) $A$ \ is a $\varphi $-positive operator in $E.$

Then an embedding 
\begin{equation*}
D^{\alpha }B_{p,\theta ,\gamma }^{l,s}\left( R^{n};E\left( A\right)
,E\right) \subset B_{p,\theta ,\gamma }^{s}\left( R^{n};E\left(
A^{1-\varkappa -\mu }\right) \right)
\end{equation*}%
is continuous and there exists a constant\ $C_{\mu }$ \ $>0$, depending only
on $\mu $, such that 
\begin{equation}
\left\Vert \eta \left( t\right) D^{\alpha }u\right\Vert _{B_{p,\theta
,\gamma }^{s}\left( R^{n};E\left( A^{1-\varkappa -\mu }\right) \right) \leq }
\tag{6.1}
\end{equation}%
\begin{equation*}
C_{\mu }\left[ h^{\mu }\left\Vert u\right\Vert _{B_{p,\theta ,\gamma
,t}^{l,s}\left( R^{n};E\left( A\right) ,E\right) }+h^{-\left( 1-\mu \right)
}\left\Vert u\right\Vert _{B_{p,\theta ,\gamma }^{s}\left( R^{n};E\right) }%
\right]
\end{equation*}%
for all $u\in B_{p,\theta ,\gamma }^{l,s}\left( R^{n};E\left( A\right)
,E\right) $ and $0<$\ $h\leq h_{0}<\infty .$

\ \ \textbf{Proof}. We have 
\begin{equation}
\left\Vert D^{\alpha }u\right\Vert _{B_{p,\theta ,\gamma }^{s}\left(
R^{n};E\left( A^{1-\varkappa -\mu }\right) \right) }=\left\Vert
A^{1-\varkappa -\mu }D^{\alpha }u\right\Vert _{B_{p,\theta ,\gamma
}^{s}\left( R^{n};E\right) }  \tag{6.2}
\end{equation}%
for all $u$ such that 
\begin{equation*}
\left\Vert D^{\alpha }u\right\Vert _{B_{p\ddot{o}\theta ,\gamma }^{s}\left(
R^{n};E\left( A^{1-\varkappa -\mu }\right) \right) }<\infty .
\end{equation*}%
On the other hand\ by using the relation $\left( 6.2\right) $ we have

\begin{equation}
A^{1-\alpha -\mu }D^{\alpha }u=F^{-1}FA^{1-\varkappa -\mu }D^{\alpha
}u=F^{-1}A^{1-\varkappa -\mu }FD^{\alpha }u=  \notag
\end{equation}%
\begin{equation}
F^{-1}A^{1-\varkappa -\mu }\left( i\xi \right) ^{\alpha }Fu=F^{-1}\left(
i\xi \right) ^{\alpha }A^{1-\varkappa -\mu }Fu.  \tag{6.3}
\end{equation}%
Hence denoting $Fu$ by $\hat{u},$ we get from the relations $\left(
6.2\right) $ and $\left( 6.3\right) $%
\begin{equation*}
\left\Vert D^{\alpha }u\right\Vert _{B_{p,\theta ,\gamma }^{s}\left(
R^{n};E\left( A^{1-\varkappa -\mu }\right) \right) }\backsim \left\Vert
F^{-1}\left( i\xi \right) ^{\alpha }A^{1-\varkappa -\mu }\hat{u}\right\Vert
_{B_{p,\theta ,\gamma }^{s}\left( R^{n};E\right) }.
\end{equation*}%
Similarly, from definition of $Y$ we have 
\begin{equation*}
\left\Vert u\right\Vert _{Y}=\left\Vert u\right\Vert _{B_{p,\theta ,\gamma
}^{s}\left( R^{n};E\left( A\right) \right) }+\sum\limits_{k=1}^{n}\left\Vert
t_{k}D_{k}^{l_{k}}u\right\Vert _{B_{p,\theta ,\gamma }^{s}\left(
R^{n};E\right) }=
\end{equation*}%
\begin{eqnarray*}
&&\left\Vert F^{-1}\hat{u}\right\Vert _{B_{p,\theta ,\gamma }^{s}\left(
R^{n};E\left( A\right) \right) }+\sum\limits_{k=1}^{n}\left\Vert t_{k}F^{-1}%
\left[ \left( i\xi _{k}\right) ^{l_{k}}\hat{u}\right] \right\Vert
_{B_{p,\theta ,\gamma }^{s}\left( R^{n};E\right) }\backsim \\
&&\left\Vert F^{-1}A\hat{u}\right\Vert _{B_{p,\theta ,\gamma }^{s}\left(
R^{n};E\right) }+\sum\limits_{k=1}^{n}\left\Vert t_{k}F^{-1}\left[ \left(
i\xi _{k}\right) ^{l_{k}}\hat{u}\right] \right\Vert _{B_{p,\theta ,\gamma
}^{s}\left( R^{n};E\right) }
\end{eqnarray*}%
for all $u\in Y.$ Thus proving the inequality\ $\left( 6.1\right) $ for some
constants $C_{\mu }$ is equivalent to proving 
\begin{equation*}
\eta \left\Vert F^{-1}\left( i\xi \right) ^{\alpha }A^{1-\varkappa -\mu }%
\hat{u}\right\Vert _{B_{p,\theta ,\gamma }^{s}\left( R^{n};E\right) }\leq
\end{equation*}%
\begin{eqnarray*}
&&C_{\mu }\left[ h^{\mu }\left( \left\Vert F^{-1}A\hat{u}\right\Vert
_{B_{p,\theta ,\gamma }^{s}\left( R^{n};E\right)
}+\sum\limits_{k=1}^{n}\left\Vert t_{k}F^{-1}\left[ \left( i\xi _{k}\right)
^{l_{k}}\hat{u}\right] \right\Vert _{B_{p,\theta ,\gamma }^{s}\left(
R^{n};E\right) }\right) \right. + \\
&&\left. h^{-\left( 1-\mu \right) }\left\Vert F^{-\shortmid }\hat{u}%
\right\Vert _{B_{p,\theta ,\gamma }^{s}\left( R^{n};E\right) })\right] .
\end{eqnarray*}%
Thus the inequality $\left( 6.1\right) $ will be followed if we prove the
following \ inequality 
\begin{equation}
\eta \left\Vert F^{-1}\left[ \left( i\xi \right) ^{\alpha }A^{1-\varkappa
-\mu }\hat{u}\right] \right\Vert _{B_{p,\theta ,\gamma }^{s}\left(
R^{n};E\right) }\leq  \tag{6.4}
\end{equation}%
\begin{equation*}
C_{\mu }\left\Vert F^{-1}\left[ h^{\mu }(A+\psi \left( t,\xi \right) \right] 
\hat{u}\right\Vert _{B_{p,\theta ,\gamma }^{s}\left( R^{n};E\right) }
\end{equation*}%
for a suitable $C_{\mu }>0$ and for all $u\in Y,$ where 
\begin{equation*}
\psi =\psi \left( t,\xi \right) =\sum\limits_{k=1}^{n}t_{k}\left\vert \xi
_{k}\right\vert ^{l_{k}})+h^{-1}.
\end{equation*}

\ Let us express the left hand side of $\left( 6.3\right) $ as follows 
\begin{equation}
\eta \left\Vert F^{-1}\left[ \left( i\xi \right) ^{\alpha }A^{1-\varkappa
-\mu }\hat{u}\right] \right\Vert _{B_{q,\theta ,\gamma }^{s}\left(
R^{n};E\right) }=  \tag{6.5}
\end{equation}%
\begin{equation*}
\eta \left\Vert F^{-1}\left( i\xi \right) ^{\alpha }A^{1-\varkappa -\mu } 
\left[ h^{\mu }(A+\psi \right] ^{-1}\left[ h^{\mu }\left( A+\psi \right) %
\right] \right\Vert _{B_{q,\theta ,\gamma }^{s}\left( R^{n};E\right) }.
\end{equation*}%
(Since\ $A$ is a positive operator in\ $E$ and $-\psi \left( t,\xi \right)
\in S\left( \varphi \right) $ so it is possible ). It is clear that the
inequality $\left( 6.4\right) $ will be followed immediately from $\left(
6.5\right) $ if we can prove that the operator-function

\begin{equation*}
\Psi _{t}=\Psi _{t,h,\mu }=\eta \left( t\right) \left( i\xi \right) ^{\alpha
}A^{1-\varkappa -\mu }\left[ h^{\mu }\left( A+\psi \right) \right] ^{-1}
\end{equation*}%
is a multiplier in $M_{p,\theta ,\gamma }^{p,\theta ,\gamma }\left( E\right)
,$ which is uniformly with respect to\ $h$ and $t.$ In order to prove that $%
\Psi _{t}\in M_{p,\theta ,\gamma }^{p,\theta ,\gamma }\left( E\right) $ it
suffices to show that there exists a constant \ $M_{\mu }>0$ with 
\begin{equation}
\left\vert \xi \right\vert ^{k}\left\Vert D^{\beta }\Psi _{t}\left( \xi
\right) \right\Vert _{L\left( E\right) }\leq C,\text{ }k=0,1,...,\left\vert
\beta \right\vert \   \tag{6.6}
\end{equation}%
for all 
\begin{equation*}
\beta =\left( \beta _{1},\beta _{2},...,\beta _{n}\right) ,\text{ }\beta
_{k}\in \left\{ 0,1\right\} ,\text{ }\xi _{k}\neq 0.
\end{equation*}%
To see this, we apply Lemma 6.1 and get a constant \ $M_{\mu }>0$ depending
only on $\mu $ such that 
\begin{equation*}
\ \left\Vert \Psi _{t}\left( \xi \right) \right\Vert _{L\left( E\right)
}\leq M_{\mu }\ \ \ \ \ \ \ \ \ \ 
\end{equation*}%
for all $\xi \in R^{n}.$ This shows that the inequality\ $\left( 7.6\right) $
is satisfied for $\beta =\left( 0,...,0\right) .$ We next consider $\left(
6.6\right) $ for $\beta =\left( \beta _{1},...\beta _{n}\right) $ \ where \ $%
\beta _{k}=1$ and $\beta _{k}=0$ for $j\neq k.$ By using the condition $%
\varkappa +\nu \left( l\right) \leq 1$ and well known inequality 
\begin{equation*}
y_{1}^{\alpha _{1}}y_{2}^{\alpha _{2}}...y_{n}^{\alpha _{n}}\leq C\left[
1+\dsum\limits_{k=1}^{n}y_{k}^{l_{k}}\right] ,\text{ }y_{k}\geq 0,
\end{equation*}%
we have 
\begin{equation*}
\left\vert \xi \right\vert \left\vert \xi _{k}\right\vert \left\Vert
D_{k}\Psi _{t}\left( \xi \right) \right\Vert _{L\left( E\right) }\leq M_{\mu
},\text{ }k=1,2...n.
\end{equation*}

Repeating the above process we obtain the estimate $\left( 7.6\right) .$
Thus the operator-function\ $\Psi _{t,h,\mu }\left( \xi \right) $\ is a
uniform collection of multiplier with respect\ to $h$ and $t$ i.e 
\begin{equation*}
\Psi _{t,h,\mu }\in \Phi _{h}\subset M_{p,\theta ,\gamma }^{p,\theta ,\gamma
}\left( E\right) .
\end{equation*}%
\ This completes the proof of the Theorem 6.1.\ It is possible to state
Theorem 6.1 in a more general setting. For this, we use the conception of
extension operator.

\ \textbf{Condition 6.1}. Let $\gamma \in A_{\nu }$ for $\nu \in \left[
1,\infty \right] $. $\ $Assume $E$ is a Banach spaces with weighted Fourier
type $\gamma $ and $\sigma \in \left[ 1,2\right] $. Suppose $A$ is a $%
\varphi $-positive operator in Banach spaces $E.$ Let a region \ $\Omega
\subset R^{n}$ be such that there exists a bounded linear extension operator$%
\ B$ from $B_{p,\theta ,\gamma }^{l,s}\left( \Omega ;E\left( A\right)
,E\right) $ to $B_{p,\theta ,\gamma }^{l,s}$\ $\left( R^{n};E\left( A\right)
,E\right) ,$ for $p$, $\theta \in \left[ 1,\infty \right] .$

\textbf{Remark 7.1}. If \ $\Omega \subset R^{n}$ is a region satisfying a
strong $l$-horn condition (see $\left[ \text{4}\right] $, \S\ 18) $E=R,$ $%
A=I $, then there exists a bounded\ linear extension operator from $%
B_{p,\theta }^{s}\left( \Omega \right) =B_{p,\theta }^{s}\left( \Omega ;%
\mathbb{C},\mathbb{C}\right) $ to 
\begin{equation*}
B_{p,\theta }^{s}\left( R^{n}\right) =B_{p,\theta }^{s}\left( R^{n};\mathbb{C%
},\mathbb{C}\right) .
\end{equation*}

Let 
\begin{equation*}
Y=B_{p,\theta ,\gamma }^{s}\left( R^{n};E\right) \text{, }Y_{0}=B_{p,\theta
,\gamma }^{l,s}\left( \Omega ;E\left( A\right) ,E\right)
\end{equation*}

\textbf{Theorem 6.2}. Suppose all conditions of the Theorem 6.1 and the
Condition 6.1 are hold. Then the embedding 
\begin{equation*}
D^{\alpha }B_{p,\theta ,\gamma }^{l,s}\left( \Omega ;E\left( A\right)
,E\right) \subset B_{q,\theta ,\gamma }^{s}\left( \Omega ;E\left(
A^{1-\varkappa -\mu }\right) \right)
\end{equation*}%
is continuous and there exists a constant $C_{\mu }$ depending only on $\mu $
such that 
\begin{equation}
\eta \left\Vert D^{\alpha }u\right\Vert _{B_{q,\theta ,\gamma }^{s}\left(
\Omega ;E\left( A^{1-\varkappa -\mu }\right) \right) }\leq  \tag{6.7}
\end{equation}%
\begin{equation*}
C_{\mu }\left[ h^{\mu }\left\Vert u\right\Vert _{B_{p,\theta ,\gamma
,t}^{l,s}\left( \Omega ;E\left( A\right) ,E\right) }+h^{-\left( 1-\mu
\right) }\left\Vert u\right\Vert _{B_{p,\theta ,\gamma }^{s}\left( \Omega
;E\right) }\right]
\end{equation*}%
for all $u\in Y_{0}$ and $0<h\leq h_{0}<\infty .$

\ \textbf{Proof}. It suffices to prove the estimate $\left( 7.7\right) .$
Let\ $P$ be a bounded linear extension operator from \ $B_{q,\theta ,\gamma
}^{s}\left( \Omega ;E\right) $ \ to $B_{q,\theta ,\gamma }^{s}\left(
R^{n};E\right) $ and\ also from $Y_{0}$ to $B_{p,\theta ,\gamma
}^{l,s}\left( R^{n};E\left( A\right) ,E\right) $. Let $P_{\Omega }$ a
restriction operator from\ $R^{n}$ to $\Omega .$ Then for any $u\in Y$\ we
have 
\begin{equation*}
\left\Vert D^{\alpha }u\right\Vert _{B_{q,\theta ,\gamma }^{s}\left( \Omega
;E\left( A^{1-\varkappa -\mu }\right) \right) }=
\end{equation*}%
\begin{equation*}
\left\Vert D^{\alpha }P_{\Omega }Pu\right\Vert _{B_{q,\theta ,\gamma
}^{s}\left( \Omega ;E\left( A^{1-\varkappa -\mu }\right) \right) }\leq
C\left\Vert D^{\alpha }Pu\right\Vert _{B_{q,\theta ,\gamma }^{s}\left(
R^{n};E\left( A^{1-\varkappa -\mu }\right) \right) }
\end{equation*}%
\begin{eqnarray*}
&\leq &C_{\mu }\left[ h^{\mu }\left\Vert Pu\right\Vert _{B_{p,\theta ,\gamma
}^{l,s}\left( R^{n};E\left( A\right) ,E\right) }+h^{-\left( 1-\mu \right)
}\left\Vert Pu\right\Vert _{B_{p,\theta ,\gamma }^{s}\left( R^{n};E\right) }%
\right] \\
&\leq &C\mu \left[ h^{\mu }\left\Vert u\right\Vert _{B_{p,\theta ,\gamma
}^{l,s}\left( \Omega ;E\left( A\right) E\right) }+h^{-\left( 1-\mu \right)
}\left\Vert u\right\Vert _{B_{p,\theta ,\gamma }^{s}\left( \Omega ;E\right) }%
\right] .
\end{eqnarray*}%
\textbf{Result 6.1}. Let all conditions of Theorem 6.2 hold. Then for all $%
u\in Y_{0}$ we have the following multiplicative estimate 
\begin{equation}
\left\Vert D^{\alpha }u\right\Vert _{B_{q,\theta ,\gamma }^{s}\left( \Omega
;E\left( A^{1-\varkappa -\mu }\right) \right) }\leq C_{\mu }\left\Vert
u\right\Vert _{B_{p,\theta ,\gamma }^{l,s}\left( \Omega ;E\left( A\right)
,E\right) }^{1-\mu }\left\Vert u\right\Vert _{B_{p,\theta ,\gamma
}^{s}\left( \Omega ;E\right) }^{\mu }.  \tag{6.8}
\end{equation}%
Indeed setting 
\begin{equation*}
h=\left\Vert u\right\Vert _{B_{p,\theta ,\gamma }^{s}\left( \Omega ;E\right)
}.\left\Vert u\right\Vert _{B_{p,\theta ,\gamma }^{l,s}\left( \Omega
;E\left( A\right) ,E\right) }^{-1}
\end{equation*}%
\ \ in $\left( 6.7\right) $ we obtain\ $\left( 6.8\right) .$\ 

\textbf{Result 6.2}. If $\ l_{1}=l_{2}=\ldots =l_{n}=m$ then we obtain the
continuity of embedding operators in the isotropic class

\begin{equation*}
\ B_{p,\theta ,\gamma }^{m,s}\left( \Omega ;E\left( A\right) E\right) .
\end{equation*}%
\ \ For $E=\mathbb{C},$ $A=I$ we obtain the embedding of weighted Besov type
spaces

\begin{equation*}
\ D^{\alpha }B_{p,\theta ,\gamma }^{l,s}\left( \Omega \right) \subset
B_{q,\theta ,\gamma }^{s}\left( \Omega \right) .
\end{equation*}

\begin{center}
\ \ \textbf{7. Application to vector-valued functions}
\end{center}

\ Let\ $s>0$ and consider the space $\left[ \text{29, \S 1.18.2}\right] $ 
\begin{equation*}
\ l_{q}^{\sigma }=\left\{ u;\text{ \ }u=\left\{ u_{i}\right\} ,\text{ }%
i=1,2,...,\infty ,\text{ }u_{i}\in C\right\}
\end{equation*}%
with the norm$\ $%
\begin{equation*}
\left\Vert u\right\Vert _{l_{q}^{\sigma }}=\left( \sum\limits_{i=1}^{\infty
}2^{iq\sigma }\left\vert u_{i}\right\vert ^{q}\right) ^{1/q}<\infty .
\end{equation*}%
Note that \ $l_{q}^{0}=l_{q}.$ Let $A$ is an infinite matrix defined in the
space\ $l_{q}$ such that 
\begin{equation*}
D\left( A\right) =l_{q}^{\sigma },\ A=\left[ \delta _{ij}2^{\sigma i}\right]
,
\end{equation*}%
\ $\ $where $\delta _{ij}=0$, when \ $i\neq j,$ \ $\delta _{ij}=1,$ when $%
i=j,$ $i,$ $j=1,2,...,\infty .$

It is clear to see that this operator $A$\ is positive in the space \ $%
l_{q}. $ Then by Theorem 7.2 we obtain the continuous embedding 
\begin{equation*}
D^{\alpha }B_{p_{1},\theta ,\gamma }^{l,s}\left( \Omega ;l_{q}^{\sigma
},l_{q}\right) \subset B_{p_{2},\theta ,\gamma }^{s}\left( \Omega
;l_{q}^{\sigma \left( 1-\varkappa -\mu \right) }\right) ,\text{ }\varkappa
=\sum\limits_{k=1}^{n}\frac{\alpha _{k}+\frac{1}{p_{1}}-\frac{1}{p_{2}}}{%
l_{k}}
\end{equation*}%
and the accociate estimate $\left( 6.7\right) ,$ where $0\leq \mu +\nu
\left( l\right) \leq 1-\varkappa .$

It should be not that the above embedding haven't been obtained with
classical method up to this time.

\begin{center}
\textbf{\ \ \ 8. }$\mathbf{B}$\textbf{-separable DOE in }$R^{n}$
\end{center}

\bigskip\ Let us consider the differential-operator equation $\left(
1.1\right) .$

\textbf{Condition 8.1. }\ Let 
\begin{equation*}
\text{(a) }K\left( \xi \right) =\sum\limits_{\left\vert \alpha \right\vert
=2l}a_{\alpha }\left( i\xi _{1}\right) ^{\alpha _{1}}\left( i\xi _{2}\right)
^{\alpha _{2}}...\left( i\xi _{n}\right) ^{\alpha _{n}}\in S\left( \varphi
\right) ;
\end{equation*}%
(b) There exists the positive constat $M_{0}$ so that%
\begin{equation*}
\text{ }\left\vert K\left( \xi \right) \right\vert \geq
M_{0}\dsum\limits_{k=1}^{n}\xi _{k}^{2l}\text{ for all }\xi \in R^{n},\text{ 
}\xi \neq 0.
\end{equation*}

\textbf{Definition 8.1}$.$ The problem $\left( 1.1\right) $ is said to be
weighted $B$-separable (or weighted $B_{p,\theta ,\gamma }^{s}\left(
R^{n};E\right) $-separable) if the problem $\left( 1.1\right) $ has a unique
solution $u\in B_{q,\theta ,\gamma }^{2l,s}\left( R^{n};E\left( A\right)
,E\right) $ for all $f\in B_{q,\theta ,\gamma }^{s}\left( R^{n};E\right) $
and 
\begin{equation*}
\left\Vert Au\right\Vert _{B_{q,\theta ,\gamma }^{s}\left( \Omega ;E\right)
}+\sum\limits_{\left\vert \alpha \right\vert =2l}\left\Vert D^{\alpha
}u\right\Vert _{B_{q,\theta ,\gamma }^{s}\left( \Omega ;E\right) }\leq
C\left\Vert f\right\Vert _{B_{q,\theta ,\gamma }^{s}\left( \Omega ;E\right)
}.
\end{equation*}

Consider the following degenerate DOE 
\begin{equation}
\ Lu=\sum\limits_{\left\vert \alpha \right\vert =2l}a_{\alpha }D^{\left[
\alpha \right] }u+Au+\sum\limits_{\left\vert \alpha \right\vert
<2l}A_{\alpha }D^{\left[ \alpha \right] }u=f  \tag{8.1}
\end{equation}%
where $A\left( x\right) $, $A_{\alpha }\left( x\right) $ are possible
unbounded operators in a Banach space $E,$ $a_{k}$ are complex-valued
functions and 
\begin{equation*}
D_{x_{k}}^{\left[ i\right] }=\left( \gamma \left( x_{k}\right) \frac{%
\partial }{\partial x_{k}}\right) ^{i},\text{ }D^{\left[ \alpha \right]
}=D_{1}^{\left[ \alpha _{1}\right] }D_{2}^{\left[ \alpha _{2}\right]
}...D_{n}^{\left[ \alpha _{n}\right] }.
\end{equation*}

\textbf{Remarke 8.1. }Under the substitution 
\begin{equation}
\tau _{k}=\int\limits_{0}^{x_{k}}\gamma ^{-1}\left( y\right) dy  \tag{8.2}
\end{equation}%
spaces $B_{p,\theta ,\gamma }^{s}\left( R^{n};E\right) ,$ $B_{p,\theta
,\gamma }^{\left[ l\right] ,s}\left( R^{n};E\left( A\right) ,E\right) $ are
mapped isomorphically onto the weighted spaces $B_{p,\theta ,\tilde{\gamma}%
}^{s}\left( R^{n};E\right) ,$ $B_{p,\theta ,\tilde{\gamma}}^{l,s}\left(
R^{n};E\left( A\right) ,E\right) $, respectively, where 
\begin{equation*}
\gamma =\prod\limits_{k=1}^{n}\gamma \left( x_{k}\right) \text{, }\tilde{%
\gamma}=\tilde{\gamma}\left( \tau \right) =\prod\limits_{k=1}^{n}\gamma
\left( x_{k}\left( \tau _{k}\right) \right) .
\end{equation*}%
Moreover, under the substitution $\left( 8.2\right) $ the degenerate problem 
$\left( 8.1\right) $ is mapped to the undegenerate problem $\left(
1.1\right) $ considered in the weighted space $B_{p,\theta ,\tilde{\gamma}%
}^{s}\left( R^{n};E\right) .$

Let

\begin{equation*}
Y=B_{q,\theta ,\gamma }^{s}\left( R^{n};E\right) \text{, }Y_{0}=B_{q,\theta
,\gamma }^{2l,s}\left( R^{n};E\left( A\right) ,E\right) .
\end{equation*}

\textbf{Theorem 8.1.}\ Suppose the following conditions hold:

(1) Condition 9.1 is hold;

(3) $s>0,$ $1\leq q$, $\theta \leq \infty ,$ $k=1,2,...,n;$

(3) $\gamma \in A_{\nu }$ for $\nu \in \left[ 1,\infty \right] $. $\ E$ is a
Banach spaces with weighted Fourier type $\gamma $ and $p\in \left[ 1,2%
\right] $;

(4) $A$ is a $\varphi $-positive operator in $E$\ and 
\begin{equation*}
\ A_{\alpha }\left( x\right) A^{-\left( 1-\left\vert \alpha \right\vert -\mu
\right) }\in L_{\infty }\left( R^{n};L\left( E\right) \right) ,\text{ }0<\mu
<1-\frac{\left\vert \alpha \right\vert }{2l}.
\end{equation*}

Then for all $\ f\in Y$ and for sufficiently large $\left\vert \lambda
\right\vert $, $\lambda \in S\left( \varphi \right) $ equation $\left(
1.1\right) $ has a unique solution $u\left( x\right) \in Y_{0}$ \ and 
\begin{equation}
\sum\limits_{\left\vert \alpha \right\vert =2l}\left\Vert D^{\alpha
}u\right\Vert _{Y}+\left\Vert Au\right\Vert _{Y}\leq C\left\Vert
f\right\Vert _{Y}.  \tag{8.3}
\end{equation}%
\ \textbf{Proof}. Firstly, we will consider leading part of the equation $%
\left( 1.1\right) $ i.e. the differential-operator equation \ 
\begin{equation}
\ \ \ \ \ \ \left( L_{0}+\lambda \right) u=\sum\limits_{\left\vert \alpha
\right\vert =2l}D^{\alpha }u+Au+\lambda u=f.  \tag{8.4}
\end{equation}%
Then we apply the Fourier transform to equation $\left( 8.4\right) $ with
respect to $x=\left( x_{1},...,x_{n}\right) $ and obtain 
\begin{equation}
\sum\limits_{\left\vert \alpha \right\vert =2l}a_{\alpha }\xi ^{\alpha }\hat{%
u}\left( \xi \right) +A_{\lambda }\hat{u}\left( \xi \right) =f^{\symbol{94}%
}\left( \xi \right) .  \tag{8.5}
\end{equation}

Since $\sum\limits_{\left\vert \alpha \right\vert =2l}a_{\alpha }\xi
^{\alpha }\geq 0$ for all\ $\xi =\left( \xi _{1},...,\xi _{n}\right) \in
R^{n}$\ therefore, $\omega =\omega \left( \lambda ,\xi \right) =\lambda
+\sum\limits_{\left\vert \alpha \right\vert =2l}a_{\alpha }\xi ^{\alpha }\in
S\left( \varphi \right) $ for all \ $\xi \in R^{n},$ i.e. operator \ $%
A+\omega $\ is invertible in $E$. Hence $\left( 8.5\right) $ implies that
the solution of equation $\left( 8.4\right) $ can be represented in the form 
\begin{equation}
u\left( x\right) =F^{-1}\left( A+\omega \right) ^{-1}f^{\symbol{94}}. 
\tag{8.6}
\end{equation}%
It is clear to see that the operator- function $\varphi _{\lambda }\left(
\xi \right) =\left[ A+\omega \right] ^{-1}$ is a multiplier in $B_{p,\theta
,\gamma }^{s}\left( R^{n};E\right) $ uniformly with respect to $\lambda .$
Actually,\ by definition of the positive operator, for all \ $\xi \in R^{n}$
and\ $\lambda \geq 0$ we get 
\begin{equation*}
\left\Vert \varphi _{\lambda }\left( \xi \right) \right\Vert _{L\left(
E\right) }=\left\Vert \left( A+\omega \right) ^{-1}\right\Vert \leq M\left(
1+\left\vert \omega \right\vert \right) ^{-1}\leq M_{0}.
\end{equation*}%
Moreover, since $D_{k}\varphi _{\lambda }\left( \xi \right) =\alpha
_{k}a_{\alpha }\xi ^{\alpha }\left( A+\omega \right) ^{-2}\xi _{k}^{-1}$
then by using the resolvent properties of positive operator $A$ we have 
\begin{equation}
\left\Vert \xi _{k}D_{k}\varphi _{\lambda }\right\Vert _{L\left( E\right)
}\leq \left\vert \alpha _{k}a_{\alpha }\right\vert \xi ^{\alpha }\left\Vert
\left( A+\omega I\right) ^{-2}\right\Vert \leq M.  \tag{8.7}
\end{equation}%
Using the estimate $\left( 8.7\right) $ we show uniform estimate 
\begin{equation}
\left\vert \xi \right\vert ^{\beta }\left\Vert D_{\xi }^{\beta }\varphi
_{\lambda }\left( \xi \right) \right\Vert _{B\left( E\right) }\leq C 
\tag{8.8}
\end{equation}%
for%
\begin{equation*}
\beta =\beta _{1},...,\beta _{n}),\text{ }\beta _{i}\in \left\{ 0,1\right\} 
\text{, }\xi =\left( \xi _{1},...,\xi _{n}\right) ,\text{ }\xi _{i}\neq 0.
\end{equation*}%
\ In a similar way we prove that the operator-functions $\varphi _{\alpha
\lambda }\left( \xi \right) =\xi ^{\alpha }\varphi _{\lambda ,t},$ $%
k=1,2,..,n$ and \ $\varphi _{0\lambda }=A\varphi _{\lambda }$ satisfiy the
estimates 
\begin{equation}
\left( 1+\left\vert \xi \right\vert \right) ^{\left\vert \beta \right\vert
}\left\Vert D_{\xi }^{\beta }\varphi _{\alpha ,\lambda }\left( \xi \right)
\right\Vert _{B\left( E\right) }\leq C,\text{ }\left( 1+\left\vert \xi
\right\vert \right) ^{\left\vert \beta \right\vert }\left\Vert D_{\xi
}^{\beta }\varphi _{0,\lambda }\left( \xi \right) \right\Vert _{B\left(
E\right) }\leq C.  \tag{8.9}
\end{equation}%
Then in view of estimates $\left( 8.8\right) $ and $\left( 8.9\right) $ we
obtain that operator-functions \ $\varphi _{\lambda },$ $\varphi _{\alpha
\lambda },$ $\varphi _{0,\lambda }$ are multipliers in\ $Y.$ By $\left(
8.9\right) $ and in view of 
\begin{equation*}
\left\Vert D^{\alpha }u\right\Vert _{Y}=\left\Vert F^{-1}\xi ^{\alpha }\hat{u%
}\right\Vert _{Y}=\left\Vert F^{-1}\xi ^{\alpha }\left( A+\omega \right)
^{-1}f^{\symbol{94}}\right\Vert _{Y},
\end{equation*}%
\begin{equation*}
\left\Vert Au\right\Vert _{Y}=\left\Vert F^{-1}A\hat{u}\right\Vert
_{Y}=\left\Vert F^{-1}\left[ A\left( A+\omega \right) ^{-1}\right] f^{%
\symbol{94}}\right\Vert _{Y}.
\end{equation*}%
\ we obtain that there exists a unique solution of equation $\left(
8.4\right) $ for all \ $f\in $\ $Y$ and the uniform estimate holds \ 
\begin{equation}
\ \ \ \ \sum\limits_{\left\vert \alpha \right\vert =2l}\left\Vert D^{\alpha
}u\right\Vert _{\ Y}+\left\Vert Au\right\Vert _{\ Y}\leq C\left\Vert
f\right\Vert _{\ Y}.  \tag{8.10}
\end{equation}%
Consider the differential operator $G_{0}$\ generated by problem $\left(
8.4\right) $, that is 
\begin{equation*}
D\left( G_{0}\right) =B_{q,\theta ,\gamma }^{2l,s}\left( R^{n};E\left(
A\right) ,E\right) ,\ G_{0}u=\sum\limits_{\left\vert \alpha \right\vert
=2l}D^{\alpha }u+Au.
\end{equation*}%
The estimate $\left( 8.10\right) $ implies that the operator $G_{0}+\lambda $
for all $\lambda \geq 0$ \ has a bounded inverse from $Y$ into $Y_{0}.$\ Let 
$G$ denote\ the differential operator\ in $Y$ generated by problem $\left(
1.1\right) .$ Namely, 
\begin{equation}
D\left( G\right) =Y_{0},\text{ }Gu=G_{0}u+L_{1}u,\text{ }\
L_{1}u=\sum\limits_{\left\vert \alpha \right\vert <2l}A_{\alpha }\left(
x\right) D^{\alpha }u.  \tag{8.11}
\end{equation}

In view of (4) condition, by virtue of Theorem 6.1, for all $u\in Y$ we have
\ \ \ \ 

\begin{equation}
\left\Vert L_{1}u\right\Vert _{Y}\leq \sum\limits_{\left\vert \alpha
\right\vert <2l}\left\Vert A_{\alpha }\left( x\right) D^{\alpha
}u\right\Vert _{Y}\leq \sum\limits_{\left\vert \alpha \right\vert
<2l}\left\Vert A^{1-\frac{\left\vert \alpha \right\vert }{2l}-\mu }D^{\alpha
}u\right\Vert _{Y}\leq  \tag{8.12}
\end{equation}%
\begin{equation*}
\leq C\left[ h^{\mu }\left( \ \sum\limits_{\left\vert \alpha \right\vert
=2l}\left\Vert D^{\alpha }u\right\Vert _{Y}+\left\Vert Au\right\Vert
_{Y}\right) +h^{-\left( 1-\mu \right) }\left\Vert u\right\Vert _{Y}\right] .
\end{equation*}%
Then from estimates $\left( 8.10\right) $ and $\left( 8.12\right) $ for $%
u\in Y_{0}$ we\ obtain

\begin{equation}
\left\Vert L_{1}u\right\Vert _{Y}\leq C\left[ h^{\mu }\left\Vert
(G_{0}+\lambda )u\right\Vert _{Y}+h^{-\left( 1-\mu \right) }\left\Vert
u\right\Vert _{Y}\right] .  \tag{8.13}
\end{equation}%
Since $\left\Vert u\right\Vert _{Y}=\frac{1}{\lambda }\left\Vert \left(
G_{0}+\lambda \right) u-G_{0}u\right\Vert _{Y}$ for all $u\in Y_{0}$\ we get 
\begin{equation}
\left\Vert u\right\Vert _{_{Y}}\leq \frac{1}{\left\vert \lambda \right\vert }%
\left[ \left\Vert \left( G_{0}+\lambda \right) u\right\Vert
_{_{Y}}+\left\Vert G_{0t}u\right\Vert _{_{Y}}\right] ,  \tag{8.14}
\end{equation}%
\begin{equation*}
\left\Vert G_{0}u\right\Vert _{_{Y}}\leq C\left[ \sum\limits_{\left\vert
\alpha \right\vert =2l}\left\Vert D^{\alpha }u\right\Vert _{B_{p,\theta
,\gamma }^{s}}+\left\Vert Au\right\Vert _{B_{p,\theta ,\gamma }^{s}}\right] .
\end{equation*}%
From estimates $\left( 8.12\right) -\left( 8.14\right) $\ for all $u\in
Y_{0} $\ we obtain 
\begin{equation}
\left\Vert L_{1}u\right\Vert _{Y}\leq Ch^{\mu }\left\Vert \left(
G_{0}+\lambda \right) u\right\Vert _{_{Y}}+C_{1}\left\vert \lambda
\right\vert ^{-1}h^{-\left( 1-\mu \right) }\left\Vert \left( G_{0}+\lambda
\right) u\right\Vert _{_{Y}}.  \tag{8.15}
\end{equation}%
Then by choosing $h$ and $\lambda $ such that $Ch^{\mu }<1,$ $%
C_{1}\left\vert \lambda \right\vert ^{-1}h^{-\left( 1-\mu \right) }<1$ from $%
\left( 9.15\right) $ we obtain the uniform estimate 
\begin{equation}
\ \ \left\Vert L_{1}\left( G_{0}+\lambda \right) ^{-1}\right\Vert _{B\left(
E\right) }<1.  \tag{8.16}
\end{equation}%
Using the relation $\left( 8.11\right) $, estimates $\left( 8.10\right) $
and $\left( 8.16\right) $ and the perturbation theory of linear operators we
obtain that the differential operator\ $G+\lambda $ is invertible from $Y$
into $Y_{0}.$ This implies the estimate $\left( 8.3\right) .$

\textbf{Result 8.1. }The Theorem 8.1 implies that the differential operator $%
G$ has a resolvent operator $\left( G+\lambda \right) ^{-1}$ for $\left\vert
\arg \lambda \right\vert \leq \varphi ,$ and the following uniform estimate
holds 
\begin{equation*}
\bigskip \sum\limits_{\left\vert \alpha \right\vert \leq 2l}\left\vert
\lambda \right\vert ^{1-\frac{\left\vert \alpha \right\vert }{2l}}\left\Vert
D^{\alpha }\left( G\mathbf{+}\lambda \right) ^{-1}\right\Vert _{B\left(
Y\right) }+\left\Vert A\left( G+\lambda \right) ^{-1}\right\Vert _{B\left(
Y\right) }\leq C.
\end{equation*}

Let $Q$ denote the operator in $B_{q,\theta }^{s}\left( R^{n},E\right) $
generated by problem $\left( 8.1\right) .$ Theorem 8.1 and Remark 8.1 imply

\textbf{Result 8.2.}\bigskip\ Let all conditions of Theorem 8.1 hold. Then
for all $\ f\in B_{q,\theta }^{s}\left( R^{n},E\right) ,$ $\lambda \in
S\left( \varphi \right) $ and for sufficiently large\ $\left\vert \lambda
\right\vert $, the equation $\left( 8.1\right) $ has a unique solution $u\in
B_{q,\theta ,\gamma }^{\left[ 2l\right] ,s}\left( R^{n};E\left( A\right)
,E\right) $ and the coercive uniform estimate holds 
\begin{equation*}
\bigskip \sum\limits_{\left\vert \alpha \right\vert \leq 2l}\left\vert
\lambda \right\vert ^{1-\frac{\left\vert \alpha \right\vert }{2l}}\left\Vert
D^{\left[ \alpha \right] }\left( Q\mathbf{+}\lambda \right) ^{-1}\right\Vert
_{B\left( B_{q,\theta ,\gamma }^{s}\left( R^{n};E\right) \right) }+
\end{equation*}

\begin{equation*}
\left\Vert A\left( Q+\lambda \right) ^{-1}\right\Vert _{B\left( B_{q,\theta
,\gamma }^{s}\left( R^{n};E\right) \right) }\leq C.
\end{equation*}%
\textbf{Remark 8.1. }The\textbf{\ }Result 8.2\textbf{\ }implies that
operator $G$\textbf{\ }is positive operator in $B_{q,\theta ,\gamma
}^{s}\left( R^{n};E\right) $. Then by virtue of $\left[ \text{29, \S 1.14.5}%
\right] $ the operator $G$ for $\varphi \in \left( \frac{\pi }{2},\pi
\right) $ is a generator of an analytic semigroup in $B_{q,\theta ,\gamma
}^{s}\left( R^{n};E\right) .$

\begin{center}
\textbf{9}. \textbf{The} \textbf{Cauchy problem for degenerate parabolic DOE 
}
\end{center}

Consider the Cauchy problem for the degenerate parabolic CDOE

\begin{equation}
\partial _{t}u+\sum\limits_{\left\vert \alpha \right\vert =2l}a_{\alpha
}D_{x}^{\alpha }u+Au=f\left( t,x\right) ,\text{ }  \tag{9.1}
\end{equation}%
\begin{equation*}
u(0,x)=0\text{, }x\in R^{n}
\end{equation*}%
in $\tilde{B}_{q,r,\gamma }^{s}\left( R_{+}^{n+1};E\right) $, where $A$ is a
linear operator in a Banach space in $E$. Let%
\begin{equation*}
F=\tilde{B}_{\mathbf{q},\mathbf{\theta },\gamma }^{s}\left(
R_{+}^{n+1};E\right) =B_{q_{1},r_{1}}^{s}\left( R_{+};F\right) \text{, }%
F_{0}=B_{q,\theta ,\gamma }^{s}\left( R^{n};E\right) \text{,}
\end{equation*}

\begin{equation*}
F_{1}=\tilde{B}_{\mathbf{q},\mathbf{\theta },\gamma }^{2l,1,s}\left(
R_{+}^{n+1};E\left( A\right) ,E\right) =B_{q_{1},r_{1}}^{1,s}\left(
R_{+};D\left( G\right) ,F\right) .
\end{equation*}%
\textbf{Theorem 9.1.}\ \textbf{\ }Assume all conditions of Theorem 8.1 hold
for $\varphi \in \left( \frac{\pi }{2},\pi \right) $ and $s>0$. Then for $%
f\in F$ the problem $\left( 9.1\right) $ has a unique solution $u\in F_{1}$
satisfying 
\begin{equation}
\left\Vert D_{t}u\right\Vert _{F}+\sum\limits_{\left\vert \alpha \right\vert
=2l}\left\Vert D^{\alpha }u\right\Vert _{F}+\left\Vert Au\right\Vert
_{F}\leq C\left\Vert f\right\Vert _{F}.  \tag{9.2}
\end{equation}

\textbf{Proof.} So, the problem $\left( 9.1\right) $\ can be express as 
\begin{equation}
\frac{du}{dt}+Gu\left( t\right) =f\left( t\right) ,\text{ }u\left( 0\right)
=0,\text{ }t\in \left( 0,\infty \right) .  \tag{9.3}
\end{equation}

The Result 9.1 implies the positivity of $G$ for $\varphi \in \left( \frac{%
\pi }{2},\pi \right) .$ Then by virtue of $\left[ \text{1, Proposition 8.10}%
\right] $ we obtain that, for $f\in F$ the Cauchy problem $\left( 9.3\right) 
$ has a unique solution $u\in F_{1}$ satisfying 
\begin{equation}
\left\Vert D_{t}u\right\Vert _{B_{q_{1},\theta _{1}}^{s}\left(
R_{+};F_{0}\right) }+\left\Vert Gu\right\Vert _{B_{q_{1},\theta
_{1}}^{s}\left( R_{+};F_{0}\right) }\leq C\left\Vert f\right\Vert
_{B_{q_{1},\theta _{1}}^{s}\left( R_{+};F_{0}\right) }.  \tag{9.4}
\end{equation}

In view of Result 8. 1 the operator $G$ is separable in $F_{0},$ therefore,
the estimate $\left( 9.4\right) $ implies $\left( 9.2\right) $.

Consider now, the Cauchy problem for the degenerate parabolic CDOE

\begin{equation}
\partial _{t}u+\sum\limits_{\left\vert \alpha \right\vert =2l}a_{\alpha
}D_{x}^{\left[ \alpha \right] }u+Au=f\left( t,x\right) ,\text{ }  \tag{9.5}
\end{equation}

\begin{equation*}
u(0,x)=0\text{, }x\in R^{n}.
\end{equation*}%
Here, $B_{p,\theta ,\gamma }^{\left[ m\right] ,s}\left( \Omega ;E\right) $
denote a $E$-valued Sobolev-Besov weighted space of functions $u$ $\in
B_{q,\theta ,\gamma }^{s}\left( \Omega ;E\right) $ that have generalized
derivatives $D_{k}^{\left[ m\right] }u\in B_{q,\theta ,\gamma }^{s}\left(
\Omega ;E\right) $ with the norm 
\begin{equation*}
\left\Vert u\right\Vert _{B_{q,\theta ,\gamma }^{\left[ m\right] ,,s}\left(
\Omega ;E\right) }=\left\Vert u\right\Vert _{B_{q,\theta ,\gamma }^{s}\left(
\Omega ;E\right) }+\sum\limits_{k=1}^{n}\left\Vert D_{k}^{\left[ m\right]
}u\right\Vert _{B_{q,\theta ,\gamma }^{s}\left( \Omega ;E\right) }<\infty .
\end{equation*}

Assume $E_{0}$ is continuoisly and densely belongs to $E.$ Here, $%
B_{q,\theta ,\gamma }^{\left[ m\right] ,s}\left( \Omega ;E_{0},E\right) $
denotes the space $B_{q,\theta ,\gamma }^{s}\left( \Omega ;E_{0}\right) \cap
B_{q,\theta ,\gamma }^{\left[ m\right] ,s}\left( \Omega ;E\right) $ with the
norm 
\begin{equation*}
\left\Vert u\right\Vert _{B_{q,\theta ,\gamma }^{\left[ m\right] ,s}\left(
\Omega ;E_{0},E\right) }=\left\Vert u\right\Vert _{B_{q,\theta ,\gamma
}^{s}\left( \Omega ;E_{0}\right) }+\sum\limits_{k=1}^{n}\left\Vert D_{k}^{%
\left[ m\right] }u\right\Vert _{B_{q,\theta ,\gamma }^{s}\left( \Omega
;E\right) }<\infty .
\end{equation*}

Let

\begin{equation*}
\Phi =\tilde{B}_{\mathbf{q},\mathbf{\theta },\gamma }^{s}\left(
R_{+}^{n+1};E\right) =B_{q_{1},r_{1}}^{s}\left( R_{+};F\right) \text{, }\Phi
_{0}=B_{q,\theta ,\gamma }^{s}\left( R^{n};E\right) \text{,}
\end{equation*}

\begin{equation*}
\Phi _{1}=\tilde{B}_{\mathbf{q},\mathbf{\theta },\gamma }^{\left[ 2l\right]
,1,s}\left( R_{+}^{n+1};E\left( A\right) ,E\right) =B_{q_{1},r_{1,\gamma
}}^{1,s}\left( R_{+};D\left( Q\right) ,F\right) .
\end{equation*}

From Theorem 8.1, Result 8.2 and Remark 8.1 we obtain the following

\textbf{Result 9.1.}\bigskip\ Assume all conditions of Theorem 8.1 hold for $%
\varphi \in \left( \frac{\pi }{2},\pi \right) $ and $s>0$. Then for $f\in F$
the equation $\left( 9.5\right) $ has a unique solution $u\in \Phi _{1}$
satisfying 
\begin{equation}
\left\Vert \partial _{t}u\right\Vert _{\Phi }+\sum\limits_{\left\vert \alpha
\right\vert =2l}\left\Vert D^{\left[ \alpha \right] }u\right\Vert _{\Phi
}+\left\Vert Au\right\Vert _{\Phi }\leq C\left\Vert f\right\Vert _{\Phi }. 
\tag{9.6}
\end{equation}

\ \ \ \ \ \textbf{Remark 9.1}. \ There are a lot of positive operators in
concrete Banach spaces. Therefore, putting concrete Banach spaces\ instead
of $E$ and concrete positive differential, pseudodifferential operators, or
finite, infinite matrices, ets.\ instead of operator $A$ on DOE $\left(
8.1\right) $ and $\left( 9.5\right) $ by virtue of Results 8.2 and 9.1 we
can obtain the maximal $B_{q,\theta ,\gamma }^{s}-$regularity properties of
different class of degenerate PDEs or system of other type equations.

\begin{center}
\textbf{10. Infinite systems of anisotropic elliptic equations }
\end{center}

Consider the following infinity systems

\begin{equation}
\left( \ L+\lambda \right) u_{m}=\sum\limits_{\left\vert \alpha \right\vert
=2l}a_{\alpha }D_{x}^{\alpha }u_{m}\left( x\right) +\left[ d_{m}\left(
x\right) +\lambda \right] u_{m}\left( x\right) +  \tag{10.1}
\end{equation}

\begin{equation*}
\sum\limits_{\left\vert \alpha \right\vert <2l}\sum\limits_{k=1}^{\infty
}d_{\alpha m}\left( x\right) D^{\alpha }u_{m}=f_{m}\left( x\right) \text{, }%
x\in R^{n}\text{, }m=1,2,...,\infty .
\end{equation*}%
Let 
\begin{equation*}
Q\left( x\right) =\left\{ d_{m}\left( x\right) \right\} ,\text{ }d_{m}>0,%
\text{ }u=\left\{ u_{m}\right\} ,\text{ }Qu=\left\{ d_{m}u_{m}\right\} ,%
\text{ }m=1,2,...\infty ,
\end{equation*}

\begin{equation*}
l_{p}\left( Q\right) =\left\{ u:u\in l_{p},\left\Vert u\right\Vert
_{l_{p}\left( Q\right) }=\left\Vert Qu\right\Vert _{l_{p}}=\left(
\sum\limits_{m=1}^{\infty }\left\vert d_{m}u_{m}\right\vert ^{p}\right) ^{%
\frac{1}{p}}<\infty \right\} ,
\end{equation*}%
\begin{equation*}
x\in R^{n}\text{, }1<p<\infty .
\end{equation*}%
Let $O$ denote the differential operator in $B_{q,\theta ,\gamma }^{s}\left(
R^{n};l_{p}\right) $ generated by problem $\left( 10.1\right) $. Let 
\begin{equation*}
B=B\left( B_{q,\theta ,\gamma }^{s}\left( R^{n};l_{p}\right) \right) .
\end{equation*}

\textbf{Condition 10. 1. }Assume $\gamma \in A_{\nu }$ for $\nu \in \left[
1,\infty \right] $ and (b) assumption of Condition 8.1. is hold. There are
positive constants $C_{1}$ and $C_{2}$ so that for $\left\{ d_{j}\left(
x\right) \right\} _{1}^{\infty }\in l_{q}$ for all $x\in R^{n}$ and some $%
x_{0}\in R^{n},$%
\begin{equation*}
C_{1}\left\vert d_{j}\left( x_{0}\right) \right\vert \leq \left\vert
d_{j}\left( x\right) \right\vert \leq C_{2}\left\vert d_{j}\left(
x_{0}\right) \right\vert .
\end{equation*}

\textbf{Theorem 10.1. }Suppose the Condition 10.1 holds. Let $a_{\alpha }\in
C_{b}\left( R^{n}\right) $, $d_{m}\in C_{b}\left( R^{n}\right) ,$ $d_{\alpha
m}\in L_{\infty }\left( R^{n}\right) $ such that 
\begin{equation*}
\text{ }\max\limits_{\alpha }\sup\limits_{m}\sum\limits_{k=1}^{\infty
}d_{\alpha m}\left( x\right) d_{k}^{-\left( 1-\frac{\left\vert \alpha
\right\vert }{2m}-\mu \right) }<M,
\end{equation*}

\begin{equation*}
\text{ for all }x\in R^{n}\text{ and }0<\mu <1-\frac{\left\vert \alpha
\right\vert }{2m}.
\end{equation*}%
\ where $p\in \left( 1,\infty \right) ,$ $q,$ $\theta \in \left[ 1,\infty %
\right] .$

Then:

(a) for all $f\left( x\right) =\left\{ f_{m}\left( x\right) \right\}
_{1}^{\infty }\in B_{q,\theta ,\gamma }^{s}\left( R^{n};l_{p}\right) ,$ for $%
\left\vert \arg \lambda \right\vert \leq \varphi $ and for sufficiently
large $\left\vert \lambda \right\vert $ problem $\left( 10.1\right) $ has a
unique solution $u=\left\{ u_{m}\left( x\right) \right\} _{1}^{\infty }$
that belongs to space $B_{q,\theta ,\gamma }^{2l,s}\left( R^{n},l_{p}\left(
Q\right) ,l_{p}\right) $ and the uniform coercive estimate holds

\begin{equation}
\sum\limits_{\left\vert \alpha \right\vert \leq 2l}\left\Vert D^{\alpha
}u\right\Vert _{B_{q,\theta ,\gamma }^{s}\left( R^{n};l_{q}\right)
}+\left\Vert Qu\right\Vert _{B_{q,\theta ,\gamma }^{s}\left(
R^{n};l_{q}\right) }\leq C\left\Vert f\right\Vert _{B_{q,\theta ,\gamma
}^{s}\left( R^{n};l_{q}\right) }.  \tag{10.2}
\end{equation}

(b) For $\left\vert \arg \lambda \right\vert \leq \varphi $ and sufficiently
large $\left\vert \lambda \right\vert $ there exists a resolvent $\left(
O+\lambda \right) ^{-1}$ of operator $O$ and 
\begin{equation}
\sum\limits_{\left\vert \alpha \right\vert \leq 2l}\left\vert \lambda
\right\vert ^{1-\frac{\left\vert \alpha \right\vert }{2l}}\left\Vert
D^{\alpha }\left( O+\lambda \right) ^{-1}\right\Vert _{B}+\left\Vert Q\left(
O+\lambda \right) ^{-1}\right\Vert _{B}\leq M.  \tag{10.3}
\end{equation}

\ \textbf{Proof.} Really, let $E=l_{q},$ $A\left( x\right) $ and $A_{\alpha
}\left( x\right) $ be infinite matrices, such that

\begin{equation*}
A=\left[ d_{m}\left( x\right) \delta _{km}\right] ,\text{ }A_{\alpha }\left(
x\right) =\left[ d_{\alpha km}\left( x\right) \right] ,\text{ }k,\text{ }%
m=1,2,...\infty .
\end{equation*}%
It is clear to see that this operator $A$ is positive in $l_{p}$. Therefore,
by virtue of Theorem 9.1 we obtain that the problem $\left( 10.1\right) $
for all $f\in B_{q,\theta ,\gamma }^{s}\left( R^{n};l_{q}\right) $, for $%
\left\vert \arg \lambda \right\vert \leq \varphi $ and sufficiently large $%
\left\vert \lambda \right\vert $ has a unique solution $u$ that belongs to
space $B_{q,\theta ,\gamma }^{2l,s}\left( R^{n};l_{p}\left( Q\right)
,l_{p}\right) $ and the estimate $\left( 10.2\right) $ hold$.$ From estimate 
$\left( 10.2\right) $ we obtain $\left( 10.3\right) .$

\begin{center}
\bigskip \textbf{11. Cauchy problem for infinite systems of parabolic
equations }
\end{center}

Consider the following infinity systems of parabolic Cauchy problem

\begin{equation}
\partial _{t}u_{m}+\sum\limits_{\left\vert \alpha \right\vert =2l}a_{\alpha
}D_{x}^{\alpha }u_{m}\left( t,x\right) +d_{m}\left( x\right) u_{m}\left(
t,x\right) +  \tag{12.1}
\end{equation}

\begin{equation*}
+\sum\limits_{\left\vert \alpha \right\vert <2l}d_{\alpha m}\left( x\right)
D^{\alpha }u_{m}\left( y,x\right) =f_{m}\left( t,x\right) ,\text{ }t\in
R_{+},\text{ }x\in R^{n},
\end{equation*}%
\begin{equation*}
u_{m}\left( 0,x\right) =0,\text{\ }m=1,2,...,\infty ,
\end{equation*}%
\begin{equation*}
F=\tilde{B}_{\mathbf{q},\mathbf{\theta },\gamma }^{s}\left(
R_{+}^{n+1};l_{p}\right) =B_{q_{1},r_{1}}^{s}\left( R_{+};l_{p}\right) \text{%
, }F_{0}=B_{q,\theta ,\gamma }^{s}\left( R^{n};l_{p}\right) \text{,}
\end{equation*}

\begin{equation*}
F_{1}=\tilde{B}_{\mathbf{q},\mathbf{\theta },\gamma }^{2l,1,s}\left(
R_{+}^{n+1};D\left( O\right) ,l_{p}\right) =B_{q_{1},r_{1}}^{1,s}\left(
R_{+};D\left( O\right) ,l_{p}\right) ,
\end{equation*}%
where $O$ is the operator in $l_{p}$ generated by problem $\left(
10.1\right) $ for $\lambda =0.$

In this section we show the following

\textbf{Theorem 11.1. }Let all conditions of Theorem 10.1 are hold$.$ Then
for $f\in F$ the Cauchy problem $\left( 11.1\right) $ has a unique solution $%
u\in F_{1}$ satisfying 
\begin{equation*}
\left\Vert D_{t}u\right\Vert _{F}+\sum\limits_{\left\vert \alpha \right\vert
=2l}\left\Vert D^{\alpha }u\right\Vert _{F}+\left\Vert Au\right\Vert
_{F}\leq C\left\Vert f\right\Vert _{F}.
\end{equation*}

\textbf{Proof. }Really, let $E=l_{q},$ $A$ and $A_{k}\left( x\right) $ be
infinite matrices, such that

\begin{equation*}
A=\left[ d_{m}\left( x\right) \delta _{km}\right] ,\text{ }A_{\alpha }\left(
x\right) =\left[ d_{\alpha km}\left( x\right) \right] ,\text{ }%
k,m=1,2,...\infty .
\end{equation*}

Then the problem $\left( 11.1\right) $ can be express in a form $\left(
9.3\right) $ where $G=O$ and 
\begin{equation*}
A=\left[ d_{m}\left( x\right) \delta _{km}\right] ,\text{ }A_{\alpha }\left(
x\right) =\left[ d_{\alpha m}\left( x\right) \right] ,\text{ }%
k,m=1,2,...\infty .
\end{equation*}%
Then by virtue of Theorem 9.1 we obtain the assertion.

\begin{center}
\ \textbf{References}
\end{center}

\bigskip\ \ \ \ \ \ \ \ \ \ \ \ \ \ \ \ \ \ \ \ \ \ \ \ \ \ \ \ \ \ \ \ \ \
\ \ \ \ \ \ \ \ \ \ \ \ \ \ \ \ \ \ \ \ \ \ \ \ \ \ \ \ \ \ \ \ \ \ \ \ \ \
\ \ \ \ \ \ \ \ \ \ \ \ \ \ \ \ \ \ \ \ \ \ 

\begin{enumerate}
\item Amann H., Linear and quasi-linear equations,1, Birkhauser, Basel 1995.

\item Amann, H., Operator-valued Fourier multipliers, vector-valued Besov
spaces, and applications, Math. Nachr. 186 (1997), 5-56.

\item Agarwal. R., O' Regan, D., Shakhmurov V. B., Separable anisotropic
differential operators in weighted abstract spaces and applications, J.
Math. Anal. Appl. 338(2008), 970-983.

\item Besov, O. V., Ilin, V. P., Nikolskii, S. M., Integral representations
of functions and embedding theorems, Nauka, Moscow, 1975.

\item Bergh J., and Lofstrom J., Interpolation spaces. An introduction,
Springer-Verlag, Berlin, 1976, Grundlehren der Mathematischen
Wissenschaften, no. 223.

\item Bourgain, J., A Hausdor -Young inequality for B-convex Banach spaces,
Paci c J. Math. 101(1982) (2), 255-262.

\item Chill, R., Fiorenza, A., Singular integral operators with
operator-valued kernels, and extrapolation of maximal regularity into
rearrangement invariant Banach function spaces. J. Evol. Equ. 14 (2014),
4-5, 795--828.

\item Diestel, J., Jarchow H., Tonge A., Absolutely summing operators,
Cambridge Univ. Press, Cambridge, 1995.

\item Denk, R., Hieber M., Pr\"{u}ss J., $R-$boundedness, Fourier
multipliers and problems of elliptic and parabolic type, Mem. Amer. Math.
Soc. 166 (2003), n.788.

\item Haller, R., Heck H., Noll A., Mikhlin's theorem for operator-valued
Fourier multipliers in $n$ variables, Math. Nachr. 244, (2002), 110-130.

\item Girardi, M., Lutz, W., Operator-valued Fourier multiplier theorems on
Besov spaces, Math. Nachr., 251, 34--51, 2003.

\item Girardi, M., Lutz, W., Operator-valued Fourier multiplier theorems on
Lp(X) and geometry of Banach spaces. J. Funct. Anal., 204(2), 320--354, 2003.

\item H\"{a}nninen, T. S., Hyt\"{o}nen, T. P., The A2 theorem and the local
oscillation decomposition for Banach space valued functions, J. Operator
Theory, 72 (2014), (1), 193--218.

\item Kurtz, D. G., Whedeen, R. L., Results on weighted norm inequalities
for multipliers, Trans. Amer. Math. Soc. 255(1979), 343--362.

\item Lions J. L and Peetre J., Sur une classe d'espaces d'interpolation,
Inst. Hautes Etudes Sci. Publ. Math., 19(1964), 5-68.

\item Meyries, M., Veraar, M., Pointwise multiplication on vector-valued
function spaces with power weights. J. Fourier Anal. Appl. 21 (2015)(1),
95--136.

\item McConnell Terry R., On Fourier Multiplier Transformations of
Banach-Valued Functions, Trans. Amer. Mat. Soc. 285, (2) (1984), 739-757.

\item Muckenhoupt, B., Hardy's inequality with weights, Studia Math., 44(1)
(1972), 31-38.

\item Pelczy\'{n}ski A., and Wojciechowski M., Molecular decompositions and
embedding theorems for vector-valued Sobolev spaces with gradient norm,
Studia Math. 107 (1993) (1), 61-100.

\item Stein, E., Singular integrals and differentiability properties of
functions, Princeton University Press (1970).

\item Shakhmurov, V. B., Theorems about of compact embedding and
applications, Doklady Akademii Nauk SSSR, 241(6), (1978), 1285-1288.

\item Shakhmurov, V. B., Imbedding theorems and their applications to
degenerate equations, Differential equations, 24 (4), (1988), 475-482.

\item Shakhmurov, V. B., Embedding operators and maximal regular
differential-operator equations in Banach-valued function spaces, Journal of
Inequalities and Applications, 4(2005), 329-345.

\item Shakhmurov, V. B., Linear and nonlinear abstract equations with
parameters, Nonlinear Anal., 73(2010), 2383-2397.

\item Shakhmurov, V. B., Coercive boundary value problems for regular
degenerate differential-operator equations, J. Math. Anal. Appl., 292 ( 2),
(2004), 605-620.

\item Shakhmurov, V. B., Embedding and maximal regular differential
operators in Sobolev-Lions spaces, Acta Mathematica Sinica, 22(5) 2006,
1493-1508.

\item Schmeisser H., Vector-valued Sobolev and Besov spaces, Seminar
analysis of the Karl-Weierstra -Institute of Mathematics 1985/86 (Berlin,
1985/86), Teubner, Leipzig, 1987, 4-44.

\item Triebel, H., Interpolation theory. Function spaces. Differential
operators, North-Holland, Amsterdam, 1978.

\item Triebel, H., Spaces of distributions with weights. Multipliers in $%
L_{p}$-spaces with weights, Math. Nachr. 78, (1977), 339-356.

\item Weis, L., Operator-valued Fourier multiplier theorems and maximal $%
L_{p}$ regularity, Math. Ann. 319, (2001), 735-75.

\item Yakubov, S and Yakubov Ya., \textquotedblright Differential-operator
equations. Ordinary and Partial \ Differential equations \textquotedblright
, Chapman and Hall /CRC, Boca Raton, 2000

\item Zimmerman, F., On vector-valued Fourier multiplier theorems, Studia
Math. 93 (3)(1989), 201-222.
\end{enumerate}

\bigskip

\begin{center}
\bigskip

\bigskip
\end{center}

\ \ \ \ \ \ \ \ \ \ \ \ \ \ \ \ \ \ \ \ \ \ \ \ 

\ 

\begin{center}
\bigskip
\end{center}

\end{document}